\documentclass[12pt]{amsart}

\usepackage{mathtools}

\usepackage{bm,cite}
\usepackage{amssymb}
\usepackage{mathrsfs}

\usepackage{graphicx}
\usepackage{subfig}
\usepackage{tabularx}
\usepackage{helvet}
\usepackage{courier}
\usepackage{type1cm}
\usepackage{color}
\usepackage{booktabs}
\usepackage{multirow}
\usepackage{geometry,calc}
\usepackage{enumerate}
\usepackage{esint}

\newtheorem{theorem}{Theorem}[section]
\newtheorem{corollary}[theorem]{Corollary}
\newtheorem{lemma}[theorem]{Lemma}
\newtheorem{proposition}[theorem]{Proposition}
\theoremstyle{definition}

\newtheorem{problem}[theorem]{Problem}
\newtheorem{assumption}[theorem]{Assumption}

\theoremstyle{remark}
\newtheorem{remark}[theorem]{Remark}

\numberwithin{equation}{section}

\newcommand{\HDGA}{$\mathrm{HDG}(k-1,k)$}
\newcommand{\HDGk}{$\mathrm{HDG}(k,k)$}

\newcommand{\HDGku}{$\text{HDG}(1,1)$}
\newcommand{\HDGkd}{$\text{HDG}(2,2)$}

\newcommand{\Tmax}{T_h^{\max}}
\newcommand{\Lu}{\Lambda_u}
\newcommand{\Lp}{\Lambda_\p}
\newcommand{\Cinv}{C_{\mathrm{inv}}}

\renewcommand{\o}{\mathfrak d}
\renewcommand{\v}{\mathfrak r}

\renewcommand{\theta}{\vartheta}
\renewcommand{\phi}{\varphi}
\renewcommand{\O}{\Omega}

\newcommand{\dt}{\partial_t}
\newcommand{\p}{\varphi}
\newcommand{\n}{\mathbf{n}}

\newcommand{\dert}{\partial_t}

\newcommand{\Au}{\mathcal{A}_u}
\newcommand{\Ap}{\mathcal{A}_{\p}}
\newcommand{\C}{\mathcal{C}}

\newcommand{\Vh}{\mathbb{V}_h}
\newcommand{\Sh}{\mathbb{S}_h}

\newcommand{\vareps}{\varepsilon}

\newcommand{\bs}[1]{\boldsymbol{#1}}

\newcommand{\fp}{\bs{\psi}}
\newcommand{\fu}{\bs{j}}
\newcommand{\hfp}{\widehat{\bs{\psi}}}
\newcommand{\hfu}{\widehat{\bs{j}}}
\newcommand{\ft}{\bs{\tau}}
\newcommand{\fv}{\bs{q}}
\newcommand{\fw}{\bs{h}}
\newcommand{\hu}{\widehat u}
\newcommand{\hp}{\widehat {\p}}
\renewcommand{\t}{\theta}
\newcommand{\ct}{\widehat\theta}
\newcommand{\hw}{\widehat w}
\newcommand{\GN}{\Gamma}

\newcommand{\tu}{\gamma_u}
\newcommand{\tp}{\gamma_\phi}

\newcommand{\Th}{{\mathcal{T}_h}}
\newcommand{\Eh}{{\mathcal{E}_h}}

\newcommand{\Hdiv}[1]{H(\mathrm{div},#1)}
\newcommand{\Poly}[1]{\mathbb{P}_{#1}}

\newcommand{\inu}{\beta}

\newcommand{\bK}{{\partial K}}

\renewcommand{\b}{\widehat e}

\newcommand{\V}{\mathbb{V}}

\newcommand{\uI}{u^I}
\newcommand{\fuI}{\fu^I}
\newcommand{\huI}{\hu^I}

\newcommand{\feta}{\boldsymbol \eta}

\newcommand{\hzeta}{\widehat \xi}
\newcommand{\fzeta}{\boldsymbol \zeta}

\newcommand{\Zeta}{\Xi}

\renewcommand{\C}[3]{\mathcal{F}(#1|#2;#3)}
\newcommand{\G}[3]{\mathcal{G}(#1|#2;#3)}
\newcommand{\B}{\mathcal{B}}

\newcommand{\chit}{\widetilde \chi}
\newcommand{\wtp}{\widetilde p}

\renewcommand{\chit}{\mu^{-1} \chi}
\renewcommand{\inu}{\mu^{-1}}

\begin{document}
	
	\allowdisplaybreaks

\addtolength{\jot}{1.5mm}
\abovedisplayskip=7pt
\belowdisplayskip=7pt

\title[HDG for a Multiple Sclerosis Model]{Hybridizable Discontinuous Galerkin Discretization of a Chemotactic Model of Multiple Sclerosis}
\author{S. Bertoluzza}
\address[Silvia Bertoluzza]{Istituto di Matematica Applicata e Tecnologie Informatiche "E. Magenes", CNR, Pavia, Italy}\email{silvia.bertoluzza@imati.cnr.it}

\author{M. Pennacchio}
\address[Micol Pennacchio]{Istituto di Matematica Applicata e Tecnologie Informatiche, "E. Magenes", CNR, Pavia, Italy}\email{micol.pennacchio@imati.cnr.it}

\thanks{Work funded by the Italian Ministry of Research (MUR) under the complementary actions to the NRRP  “D34 Health: Digital Driven Diagnostics, prognostics, and therapeutics for sustainable Health care” Grant  \#PNC0000001, CUP B53C22006100001. MP has been partially supported by ICSC-Centro Nazionale di Ricerca in High Performance Computing, Big Data, and Quantum Computing funded by the European Union--NextGenerationEU. The authors are members of the Istituto Nazionale di Alta Matematica - Gruppo Nazionale di Calcolo Scientifico.}

\keywords{hybridizable discontinuous Galerkin, Keller--Segel model, error estimates, chemotaxis, multiple sclerosis, pattern formation}
\subjclass[2020]{65M60, 65M12, 65M15, 35K57, 92C17}


\maketitle

\begin{abstract}
	We propose and analyze a hybridizable discontinuous Galerkin (HDG) method for
	a chemotaxis model of multiple sclerosis, in which a Keller--Segel type system
	for activated macrophages and a chemoattractant is coupled with an ordinary
	differential equation for the density of destroyed oligodendrocytes. All the
	unknowns, including the total cellular flux, are approximated with the same
	polynomial degree $k$. We prove existence and uniqueness of the semi-discrete
	solution, and convergence of order $k+1$ in $L^2$ for all the variables, fluxes
	included. Numerical experiments confirm the theoretical rates and show that
	the method reproduces the patterns known for this model, also on a curved,
	non-convex domain.
\end{abstract}

\section{Introduction}

Mathematical models and their numerical simulation are increasingly used to
 investigate the mechanisms underlying diseases of the central nervous
system, 
from the spreading of misfolded proteins in neurodegenerative diseases such as
Parkinson's disease \cite{Antonietti2024heterodimer,Corti2023FK} to fluid flow
and solute transport in the brain \cite{Corti2023poro,MardalRognes2022}.
A common requirement of
these simulations is a discretization that is accurate on the complex
geometries arising from the anatomy of the brain.

Here we focus on multiple sclerosis (MS), a chronic inflammatory disease of the central nervous
system, in which the myelin sheath of the axons is damaged and the
oligodendrocytes, the cells that produce it, are destroyed. In one of the proposed
mechanisms for the early stage of the lesions, activated macrophages attack the
oligodendrocytes, and the damaged tissue in turn releases chemical signals that
recruit further macrophages. A mathematical model of this mechanism was proposed
in \cite{CALVEZ,CalvezKhonsari2008} and further developed in
\cite{lombardo2017demyelination}: the density $u$ of activated macrophages and the
concentration $\varphi$ of a chemoattractant obey a Keller--Segel type system
\cite{KS}, coupled with an ordinary differential equation for the density
$\mathfrak d$ of destroyed oligodendrocytes, which acts as a source of
chemoattractant. The well-posedness of the model, in the more general form
considered here, was proved in \cite{desvillettes2022global}. Above a critical
value of the chemotactic sensitivity the homogeneous equilibrium undergoes a
Turing instability \cite{turing1990chemical}, and the resulting stationary patterns
have been related to the different forms of demyelinated plaques
observed in the disease \cite{lombardo2017demyelination,Bisi_pattern}.

To our knowledge, the numerical simulations available for this model are
restricted to square domains, or to radially symmetric configurations, with
no-flux or periodic boundary conditions, and they mainly rely on finite
difference or spectral discretizations; see, e.g.,
\cite{CalvezKhonsari2008,lombardo2017demyelination,desvillettes2022global,Bisi_pattern}.
On the other hand, the shape of the region in which a lesion develops does
matter. Concentric lesions are observed to turn into bands where the tissue is
narrower, as in the brainstem \cite{CALVEZ}; the inclusion of the cerebral
ventricles, in order to obtain a more realistic two dimensional domain, is
indicated in \cite{lombardo2017demyelination} as the main improvement needed by
the available simulations; and, more generally, the patterns selected by a
Turing instability depend on 
 the size and shape of the domain \cite{murray_mathematical_biology,Arcuri_Murray_1986},  and on the boundary conditions \cite{Arcuri_Murray_1986,Woolley_2022}. 
This calls for a discretization that handles
domains of general shape and general boundary conditions, and whose accuracy is
established rigorously.

Hybridizable discontinuous Galerkin (HDG) methods \cite{cockburn_unified_2009}
are well suited to this purpose. They are defined on unstructured meshes; the
fluxes are approximated together with the densities; boundary conditions of
Dirichlet, Neumann and Robin type, and not only no-flux or periodic ones, are
imposed through the numerical traces and the numerical fluxes, in the same way
as the transmission conditions between the elements; and, by static
condensation, the only globally coupled unknowns are the traces on the skeleton
of the mesh.

HDG methods for the Keller--Segel system were analyzed in \cite{BBPP}.  
There, as in the
method of \cite{CockburnDong2009} for linear convection--diffusion--reaction
problems, the total cellular flux, which includes the chemotactic drift, is taken
as an independent unknown.
 This allows
for an equal-order method, HDG$(k,k)$, converging with order $k+1$ in all the
variables, whose analysis relies on an energy argument only and does not require
the domain to be convex \cite[Remark~4.7]{BBPP}.
In this paper we extend the method HDG$(k,k)$ to the MS model. 

This is, to our knowledge, the first discretization of the model of
	\cite{lombardo2017demyelination,desvillettes2022global} for which a convergence
	analysis is available.
With respect to the
Keller--Segel case \cite{BBPP}, the new ingredients are the following. The
equation for $\mathfrak d$ contains no spatial derivatives: its discrete
counterpart is local to each element, it requires neither a trace nor a numerical
flux, and it leaves the structure of the condensed global system unchanged; in the
analysis, it gives rise to an additional coupling, through the damage rate $g(u)$, which is
handled by a separate bound for $\mathfrak d - \mathfrak d_h$ (Lemma~\ref{lem:boundd}).
Moreover, the chemotactic sensitivity $\chi f(u)$ and the reaction term depend
nonlinearly on the cell density, whereas in \cite{BBPP} the sensitivity is
constant.
Finally, as the discrete solution is in general not non-negative, the arguments
used for the continuous problem in \cite{desvillettes2022global}, which rely on
positivity, are not available, and all our results are proved without any sign
information on the discrete solution. On the other hand, the projection estimates
for the cellular unknowns and the bound for the chemical unknowns in terms of the
error on the densities are taken, with minor modifications, from \cite{BBPP}.

We
prove existence and uniqueness of the semi-discrete solution by a purely finite
dimensional argument and, when the chemotactic sensitivity is bounded, as in the
model of \cite{lombardo2017demyelination}, an $L^2$ stability bound which is
uniform in the mesh size (Theorem~\ref{thm:existence}). We then prove an a priori
error estimate (Theorem~\ref{thm:error}) which, for sufficiently smooth solutions,
gives the optimal order $k+1$ for the three densities and for both fluxes.

The numerical tests confirm these rates for $k=1,2$. We then show that the method
recovers the stripe, hexagon and spot patterns computed in \cite{Bisi_pattern} on
a square, and that the same simulations run, without any modification, on a
curved, non-convex domain, where the patterns keep their wavelength and adapt
their arrangement to the geometry. We stress that our aim is not a study of
pattern selection for the MS model, which is beyond the scope of this work, but to
provide a discretization with which such a study can be carried out on domains of
general shape.

The paper is organized as follows. In Section~\ref{sec:model} we recall the model
and the available results on its well-posedness. In Section~\ref{sec:HDG_MS} we
introduce the HDG discretization and state the existence and stability result for
the discrete solution, whose proof is given in Section~\ref{sec:existence}.
Section~\ref{sec:error} is devoted to the error analysis. In
Section~\ref{sec:numerical} we present the numerical tests, and some conclusions
are drawn in Section~\ref{sec:conclusions}.

\section{The mathematical model}\label{sec:model}
 In this paper we consider numerical simulation of the initial stages of the demyelination process through oligodendrocyte destruction driven by activated macrophages. We will consider, for such a process, the model for MS proposed by \cite{lombardo2017demyelination}, in the more general version considered by \cite{desvillettes2022global}:
\begin{gather}\label{MS1}
\dt u - \nabla \cdot \nabla u + \nabla \cdot (\chi \nabla \p f(u))  - u(1 - u^{a-1}) = 0,\qquad  \text{ in }[0,T] \times \O\\ \label{MS2}
\dt \o - g(u)(1-\o) = 0,\qquad  \text{ in }[0,T] \times \O \\[3mm] \label{MS3}
\dt \phi - \mu \nabla \cdot \nabla \p+ \p - (\beta u  + \delta \o ) = 0,\qquad  \text{ in }[0,T] \times \O\\[3mm] \label{MS4}
\nabla u \cdot \n = 0, \qquad \nabla \phi \cdot \n = 0 \qquad \text{ on } [0,T] \times \GN\\ \label{MS5}
u(0,\cdot) = u_0, \qquad \p(0,\cdot) = \p_0, \qquad \o(0,\cdot) = \o_0.
\end{gather}

Here $u$ is the density of immune cells, $\p$ the density of the chemoattractant and $\o$ the density of destroyed oligodendrocytes, and the initial data $u_0$, $\p_0$ and $\o_0$ are non negative. Following \cite{desvillettes2022global} we make the following assumptions on the parameters of the equation: $a > 1$, $\chi, \mu, \beta, \delta> 0$. Moreover the non linearities satisfy the following assumptions:

\newcommand{\expb}{\varsigma}
\begin{description}
\item[F]\label{assmpF} \(f \in C^2([0,\infty))\) with $f(0) = 0$ and  there exist $\expb, \ell \geq 0$  such that
\[
| f(y) |\lesssim y^\expb, \qquad | f'(y) | \lesssim y^\ell.
\]
\item[G] \(g \in C^1([0,\infty))\) and $g(y) \geq 0$.
\end{description}
 Under this assumption the following existence and uniqueness results hold:
\begin{theorem}[Desvillettes et al. \cite{desvillettes2022global}] \label{existunique}
Let $\O$ be of class $C^{2+\alpha}$, for some $\alpha > 0$, and assume that 
$\max\{ \expb-1,\expb/2\} < a / 4$. 
Let $u_0 \in W^{2-2/\wtp,\wtp}(\O)$, $\p_0 \in W^{2-2/\wtp,\wtp}(\O)$ and $\o_0 \in L^\infty(\O)$, for some $\wtp > \max\{4,a\}$. Then \eqref{MS1}--\eqref{MS5} admits a unique non negative strong solution $(u,\p,\o)$ with $\dt u$, $\dt \p$, $\dt\o$, $\Delta u$, $\Delta  \p$, $\nabla \cdot (f(u)\nabla \p)$ in $L^{\wtp}([0,T]\times\O)$ for all $T > 0$, and satisfying
\[
\| u \|_{L^\infty(0,T;L^\infty(\O))} + \| \p \|_{L^\infty(0,T;W^{1,\infty}(\O))} + \| \o \|_{L^\infty(0,T;L^\infty(\O))} \lesssim 1,
\]
the implicit constant in the above bound depending on the initial data. Moreover, if the initial data are smooth and satisfy compatibility conditions at the boundary $\nabla u_0\cdot \n = \nabla \phi_0\cdot \n =0$ on $\partial\Omega$, then the solution is classical.
\end{theorem}
Moreover, if $\chi$ is sufficiently small, the constant state is exponentially asymptotically stable. More precisely, we have the following theorem.
\begin{theorem}[Desvillettes et al. \cite{desvillettes2022global}] 
\label{thm:stability}
	Under the assumptions of Theorem \ref{existunique}, assume also that $f(1) > 0$, and that $\chi < \chi_\text{subcrit} = 4\sqrt{\mu(a-1)}/(\beta f(1))$. Then there exists an $\vareps > 0$ such that,  letting
	\[
	\bar u = 1, \qquad \bar \p = \beta + \delta, \qquad \bar \o = 1, 
	\]
	if   the initial data satisfy $u_0 \in W^{2-2/\wtp,\wtp}(\O)$, $\p_0 \in W^{2-2/\wtp,\wtp}(\O)$ and $\o_0 \in L^\infty(\O)$, and if 
\[
\| u_0-\bar u\|_{L^2(\O)} + \| \p_0 - \bar \p \|_{L^2(\O)} + \| \o_0 - \bar\o \|_{L^2(\O)}
< \varepsilon
\]
then for some $\rho > 0$ it holds that
\[
\| u(t) -\bar u\|_{L^\infty(\O)} + \| \p(t) - \bar \p \|_{W^{1,\infty}(\O)} + \| \o(t) - \bar\o \|_{L^\infty(\O)} \lesssim e^{-\rho t} \qquad \text{ for all } t\geq  0.
\]
\end{theorem}

\section{The discretization}\label{sec:HDG_MS}
 
 \subsection{Notation} In the following we let $\O \subset \mathbb{R}^2$ denote a bounded connected domain, and we let $\GN = \partial\O$ denote its boundary, $\n$ is the outer unit normal to $\GN$. We let $\Th$ denote a quasi-uniform, shape regular triangular mesh of the domain $\Omega$ and we let $\Eh$ denote  the set of edges of $\Th$. By abuse of notation, we will also use the notation $\Eh$ to indicate the skeleton of the decomposition  ($\Eh = \cup_K \partial K$). Which of the two definitions applies will be clear from the context. For $K \in \Th$ we let $h_K$ denote the diameter of $K$, and we let $h = \max_{K\in\Th} h_K$ denote the mesh size.
 For $u \in L^2(\O)$ and $\fu \in (L^2(\O))^2$ we let $u^K = u|_K$ and $\fu^K = \fu|_K$. 
 
 \
 
 We let
 \begin{gather}
 	H^1(\Th) = \{u \in L^2(\O)\ \text{such that}\ \nabla u^K \in L^2(K)\ \text{ for all } K \in \Th \},\\
 	\Hdiv{\Th} = \{\fu \in( L^2(\O))^2 \ \text{such that}\ \nabla \cdot \fu^K \in L^2(K)\ \text{ for all } K \in \Th \}.
 \end{gather} 
 
 Using the notation $(\,\cdot\,)_{K\in\Th}$ to indicate vectors with index in $\Th$, by abuse of notation, for $u \in H^1(\Th)$, $\fu \in \Hdiv{\Th}$ we let
 \begin{gather}
 	\nabla u = (\nabla u^K)_{K\in\Th} \in \prod_K (L^2(K))^2,\\
 	\nabla \cdot \fu = (\nabla \cdot \fu^K)_{K\in\Th}  \in \prod_K L^2(K), \\ \fu\cdot \n = (\fu^K \cdot \n^K)_{K\in\Th}  \in \prod_K H ^{-1/2}(\partial K)\end{gather} 
 respectively denote the piecewise gradient of $u$ and the piecewise divergence and the piecewise normal trace of $\fu$.
 We will use the following notation for the $L^2$ scalar products
 \[
 (u,v)_\O = \int_\Omega u v, \qquad (u^K,v^K)_K = \int_K u^K v^K.
 \]
 We also introduce the following compact notation for  broken scalar products: for $A, B$  in $\prod_K L^2(K)$ and $C, D$ in $\prod_K L^2(\partial K)$, 
 \[ 
 (A , B)_{\Th} = \sum_{K} \int_K A^K B^K, \qquad (C, D)_{\Eh} = \sum_K \int_{\partial K} C^K D^K
 \]
 with corresponding broken norms
 \[
 \| A \|_{\Th}^2 = (A , A)_{\Th}, \qquad \| C \|_{\Eh}^2 = (C , C)_{\Eh}.
 \]

 We observe that, with our notation, for all $u \in H^1(\Th)$, $\fu \in \Hdiv{\Th}$ the result of an element by element integration by parts can be written in compact form as
 \[
 (u,\nabla\cdot \fu)_{\Th} = -(\nabla u,\fu)_{\Th} + (u ,\fu\cdot \n)_{\Eh}.
 \]

\subsection{The discrete problem}
 In order to discretize Problem \eqref{MS1}--\eqref{MS5}, we start by considering an equivalent mixed formulation, where, for the sake of generality, we also allow for non homogeneous forcing terms in the equations corresponding  to \eqref{MS1} and \eqref{MS3} (this will allow us to test our approach on arbitrary manufactured solutions). To this aim we introduce additional flux unknowns $\fu$ and $\fp$, implicitly defined as
 \begin{equation}\label{eq:def_fluxes}
 	\fu + \nabla u + \mu^{-1} \chi f(u) \fp = 0, \qquad \mu^{-1} \fp + \nabla \p = 0,
 \end{equation}
  and we rewrite the problem as follows.
\begin{gather}\label{eq:mixed_form}
\dt u + \nabla \cdot \fu + u(u^{a-1} - 1) = f_u, \qquad \dt \p + \nabla \cdot \fp - \delta \o + \p - \beta u = f_\p\\
{\dt \o - g(u)(1-\o)} = 0.
\end{gather}
We easily see that substituting the definition of $\fu$ and $\fp$ resulting from \eqref{eq:def_fluxes} in \eqref{eq:mixed_form} yields \eqref{MS1}--\eqref{MS3} with right hand sides $f_u$, $f_\p$.
We  let 
\[
	\Poly{k}(\Th) = \{ \phi \in L^2(\O): \ \phi|_K \in \Poly{k}\ \text{for all } K \in \Th\}\]
	denote the space of discontinuous piecewise polynomials of degree at most $k$, and we let
\[	\Poly{k}(\Eh) = \{ \phi \in L^2(\Eh): \ \phi|_e \in \Poly{k}\ \text{for all } e \in \Eh\}\]
denote the space of single-valued discontinuous piecewise polynomials of degree at most $k$ on the skeleton.

\

We introduce the traces $\hu \in H^{1/2}(\Eh)$ of $u$ and $\hp \in H^{1/2}(\Eh) = H^1(\O)|_\Eh$ of $\p$ on the skeleton $\Eh$ as additional independent unknowns, which we denote by $\hu$ and $\hp$, and we proceed as usual, by introducing discontinuous test functions $w \in \Poly{k}(\Th)$, $\fv \in (\Poly{k}(\Th))^2$, $\t \in \Poly{k}(\Th)$, $\ft \in (\Poly k(\Th))^2$, $\v \in \Poly{k}(\Th)$, defined on the  mesh $\Th$. We then multiply each equation by the corresponding test function, integrate over $\O$, and, after integration by parts in each element $K$, in the element boundary terms replace the traces of $u$ and $\p$ and the normal fluxes $\fu\cdot\n$ and $\fp\cdot \n$  with numerical traces $\hu$ and $\hp$ and numerical fluxes $\hfu\cdot\n$ and $\hfp\cdot\n$, defined on $\bK$ as
\[
\hfu^K\cdot\n^K = \fu^K\cdot\n^K + \tu (u^K - \hu), \qquad \hfp^K\cdot\n^K = \fp^K\cdot\n^K + \tp(\p^K - \hp),
\]
where $\tu$ and $\tp$ are positive stabilization parameters which, as in the equal-order method of \cite{BBPP}, we take independent of the mesh size,
\[
\tu \simeq 1, \qquad \tp \simeq 1.
\]
Proceeding as described in \cite{BBPP}, we obtain the following identities \begin{gather*}
	{	(\dert u,w)_{\Omega} +  (\hfu \cdot \n, w)_{\Eh} - 
		(\fu,\nabla w)_{\Th} + 
			\underbrace{(u(u^{a-1}-1),w)_\O}_{\B(u,w)}
		= 	(f_u,w)_\Omega ,}
	\\ 
	 \dt(\o,\v) 
	+ \underbrace{(g(u)(\o-1),\v)}_{\G {u} \o \v} = 0,\\ 
	{	(\fu ,\fv)_{\Omega}
		- (u ,
		\nabla \cdot \fv)_{\Th} +  (\hu , \fv\cdot \n)_{\Eh}
		+ \underbrace{(\mu^{-1} \chi f(u) \fp ,\fv)_{\Omega}}_{\C{\fp}{u}\fv}	=0,}
	\\[3mm]
	\dert (	{\p,\t)_\Omega +	( \hfp\cdot\n, \t)_{\Eh}  -	(  \fp,\nabla\t)_{\Th}   + (\p,\t)_\Omega- (\beta u {+ \delta \o},\t)_\Omega
		= 	(f_\p,\t)_\Omega,} \\ 
	  (\inu \fp,\ft)_{\Omega}
	- (\p ,
	\nabla \cdot \ft)_{\Th} + (\hp , \ft\cdot \n)_{\Eh}
	= 0,
\end{gather*}
holding for all $w, \t  \in  H^1(\Th)$, $\fv, \ft \in \Hdiv\Th$, as well as 
\begin{equation}\label{eq:interface}
	 ( \hfu  \cdot \n, \hw)_{\Eh}  = 0,\qquad
	 -  ( \hfp \cdot \n,\ct)_{\Eh} = 0,
\end{equation}
which, provided $u$ and $\p$ are sufficiently smooth, holds for all $\hw, \ct \in L^2(\Eh)$. 	We point out that the homogeneous Neumann boundary conditions are weakly encoded in the two equations in \eqref{eq:interface}. Indeed the second equation implies $\fp\cdot \n = -\mu \nabla \p\cdot \n = 0$ on $\partial \Omega$. This in turn implies that, on $\partial\Omega$, $\fu \cdot \n = - \nabla u \cdot \n$ and, consequently, the first equation in \eqref{eq:interface} implies the homogeneous Neumann condition for $u$.  We underline that, while for the sake of simplicity we focus here on homogeneous Neumann boundary conditions, the above formulation can easily include non homogeneous Neumann and Dirichlet boundary conditions (the latter imposed directly on $\hu$ and $\hp$), as well as non homogeneous Robin boundary conditions.

As far as the discrete spaces are concerned, we approximate all the unknowns, including the cellular flux $\fu$ and the trace $\hu$, by discontinuous piecewise polynomials of the same degree $k \geq 1$. In the notation of \cite{BBPP}, where the first index denotes the polynomial degree used for the cellular flux and the second the degree used for all the other unknowns, this is the equal-order method \HDGk. We recall that in \cite{BBPP} the variant \HDGA, in which the cellular flux is approximated by polynomials of degree $k-1$, is also analyzed: for the Keller--Segel system, that variant only yields order $k$ convergence for the cellular flux, while, thanks to the definition \eqref{eq:def_fluxes} of the flux $\fu$ (which includes the chemotactic drift), the equal-order method yields order $k+1$ convergence for all the variables. 
We will then consider the following discrete problem.
\begin{problem}\label{pb:discrete}
Find $u_h, \p_h, \o_h \in \Poly k (\Th)$, $\fu_h \in (\Poly{k} (\Th))^2$, $\fp_h \in (\Poly{k}(\Th))^2$, $\hu_h, \hp_h \in \Poly {k} (\Eh)$ such that for all
$w, \t, \v \in \Poly k (\Th)$, $\fv \in (\Poly{k} (\Th))^2$, $\ft \in (\Poly{k}(\Th))^2$, $\hw, \ct \in \Poly {k} (\Eh)$ it holds that
\begin{gather}	\label{MS-HDG1}
{	(\dert u_h,w)_{\Omega} +  (\hfu_h \cdot \n, w)_{\Eh} - 
(\fu_h ,\nabla w)_{\Th} + 
\underbrace{(u_h(u_h^{a-1}-1),w)_\O}_{\B(u_h,w)}
= 	(f_u,w)_\Omega ,}
\\ 
\label{MS-HDG7} \dt(\o_h,\v)_\O
+ \underbrace{(g(u_h)(\o_h-1),\v)_\O}_{\G {u_h} {\o_h} \v} = 0,\\ 
{	\label{MS-HDG3}	(\fu_h ,\fv)_{\Omega}
- (u_h ,
\nabla \cdot \fv)_{\Th} +  (\hu_h , \fv\cdot \n)_{\Eh}
+ \underbrace{(\chit f(u_h) \fp_h ,\fv)_{\Omega}}_{\C{\fp_h}{u_h}\fv}	=0,}\\ 
\label{MS-HDG6} - ( \hfu_h  \cdot \n, \hw)_{\Eh}  = 0,
\\[8mm]
\label{MS-HDG2}	\dert (	{\p_h,\t)_\Omega +	( \hfp_h\cdot\n, \t)_{\Eh}  -	(  \fp_h,\nabla\t)_{\Th}   + ( \p_h,\t)_\Omega- (\beta u_h {+ \delta \o_h},\t)_\Omega
= 	(f_\p,\t)_\Omega,} \\ 
\label{MS-HDG4}		  (\inu \fp_h,\ft)_{\Omega}
- (\p_h ,
\nabla \cdot \ft)_{\Th} + (\hp_h , \ft\cdot \n)_{\Eh}
= 0,\\
\label{MS-HDG5}-  ( \hfp_h \cdot \n,\ct)_{\Eh} = 0,
\end{gather}
\end{problem}
As the solution to such a problem will not necessarily be non-negative and since $f$ and $g$ are only given on $[0,\infty)$, and $u^{a-1}$ is not defined for $u<0$ when $a \notin \mathbb N$, throughout the paper we will extend these non linearities as follows: for $y \in \mathbb R$ we set
\begin{equation}\label{eq:extension}
	f(y) := f(|y|), \qquad g(y) := g(|y|), \qquad b(y) := y|y|^{a-1} - y.
\end{equation}
Thanks to the finite dimensionality of the spaces involved, existence of the solution of the discrete problem can be proven by classical arguments. More precisely, letting $u_h(0), \p_h(0), \o_h(0) \in \Poly{k}(\Th)$ denote arbitrary discrete initial data (e.g., the $L^2(\O)$ orthogonal projections of $u_0$, $\p_0$ and $\o_0$ onto $\Poly{k}(\Th)$), we have the following theorem, the proof of which we postpone to {Section} \ref{sec:existence}.
\begin{theorem}[Existence and uniqueness of the discrete solution]\label{thm:existence}
	Let Assumptions {\bf F} and {\bf G} hold, with $f$, $g$ and $b$ extended to $\mathbb R$ according to \eqref{eq:extension}, and let $f_u, f_\p \in C^0([0,T];L^2(\O))$. Then, for every $h$ and for every choice of the discrete initial data $u_h(0), \p_h(0), \o_h(0) \in \Poly{k}(\Th)$, the following holds.
	\begin{enumerate}[(i)]
		\item There exist $\Tmax \in (0,T]$ and a unique solution $(u_h,\hu_h,\fu_h,\o_h,\p_h,\hp_h,\fp_h)$ of Problem \eqref{MS-HDG1}--\eqref{MS-HDG5} on $[0,\Tmax)$, with $u_h$, $\o_h$, $\p_h$, $\fp_h$ and $\hp_h$ in $C^1([0,\Tmax))$ and $\fu_h$, $\hu_h$ in $C^0([0,\Tmax))$ (with values in the respective discrete spaces). Moreover, if
		\[
		\sup_{t < \Tmax} \| u_h(t) \|_\O < \infty,
		\]
		then $\Tmax = T$ and the solution extends to a solution on the closed interval $[0,T]$, with the same regularity.
		\item For all $t \in [0,\Tmax)$ it holds that
		\begin{gather}\label{eq:boundo_stab}
			\| \o_h(t) - 1 \|_\O \leq \| \o_h(0) - 1 \|_\O,\\ \label{eq:boundp_stab}
			\begin{multlined}
				\| \p_h(t) \|_\O^2 + 2 \int_0^t \Big( \inu \| \fp_h \|_\O^2 + \tp \| \p_h - \hp_h \|_\Eh^2 \Big) \\
				\leq \| \p_h(0) \|_\O^2 + \frac 32 \int_0^t \Big( \beta^2 \| u_h \|_\O^2 + \delta^2 \| \o_h \|_\O^2 + \| f_\p \|_\O^2 \Big).
			\end{multlined}
		\end{gather}
		\item If, in addition, $f$ is bounded (that is, if $\expb = 0$ in Assumption {\bf F}), then $\Tmax = T$ and
		\begin{multline}\label{eq:stability}
			\sup_{t \leq T} \Big( \| u_h(t) \|_\O^2 + \| \p_h(t) \|_\O^2 + \| \o_h(t) \|_\O^2 \Big) \\ + \int_0^T \Big( \| \fu_h \|_\O^2 + \tu \| u_h - \hu_h \|_\Eh^2 + \inu \| \fp_h \|_\O^2 + \tp \| \p_h - \hp_h \|_\Eh^2 \Big) \leq C,
		\end{multline}
		where $C$ depends on $T$, $\chi$, $\mu$, $\beta$, $\delta$, $\| f \|_{L^\infty(0,\infty)}$, $|\O|$, on $\| u_h(0) \|_\O$, $\| \p_h(0) \|_\O$, $\| \o_h(0) \|_\O$ and on $\| f_u \|_{L^2(0,T;L^2(\O))}$, $\| f_\p \|_{L^2(0,T;L^2(\O))}$, but not on $h$.
	\end{enumerate}
\end{theorem}

\begin{remark}
	The discrete initial data can be arbitrary, and, in particular, need not be non negative. Statement (iii) covers the model of \cite{lombardo2017demyelination} considered in Section \ref{sec:numerical}, for which $f(u) = u/(1+u)$ is bounded, and yields, for that model, an $L^2$ stability bound which is uniform in $h$. For unbounded $f$ (that is, $\expb > 0$), the discrete energy estimate does not close uniformly in $h$, as the continuous argument of \cite{desvillettes2022global} relies on parabolic regularity and duality estimates which have no discrete counterpart; existence on the whole interval $[0,T]$ will then be obtained, for $h$ small enough, as a byproduct of the error analysis (see Corollary \ref{cor:global} below).
\end{remark}

\renewcommand{\S}{\mathbb{S}}

\newcommand{\e}{{\mathfrak e}}
\newcommand{\hTh}{\widehat T_h}
\newcommand{\upper}{\widetilde{K}(\o,u)}

\section{Error estimate}\label{sec:error}
We start by stating the smoothness assumptions on the solution of the continuous problem under which our error bounds will be proven. In the following, $\fp = -\mu\nabla\p$ and $\fu = -\nabla u - \chit f(u) \fp$ denote the fluxes associated with the continuous solution, as defined in \eqref{eq:def_fluxes}.
\begin{assumption}[Smoothness of the continuous solution]\label{ass:smoothness}
	The solution $(u,\o,\p)$ of \eqref{MS1}--\eqref{MS5} and the associated fluxes $\fu$ and $\fp$ satisfy
	\begin{gather}\label{eq:smooth-u}
		u \in H^1(0,T;H^s(\O)), \qquad \fu \in H^1(0,T;H^\sigma(\O)^2), \\ \label{eq:smooth-p}
		\p \in H^1(0,T;H^r(\O)), \qquad \o \in H^1(0,T;H^q(\O)),
	\end{gather}
	for some regularity indices $s$, $\sigma$, $r$, $q$ with
	\begin{equation}\label{eq:smooth-indices}
		1 < s,\sigma,q, \qquad 2 < r.
	\end{equation}
\end{assumption}
No upper bound is imposed on the regularity indices. As the approximation order of $\Poly{k}(\Th)$ is limited, only the capped indices
\begin{equation}\label{eq:capped-indices}
	\widetilde s = \min\{s,k+1\}, \qquad \widetilde \sigma = \min\{\sigma,k+1\}, \qquad \widetilde q = \min\{q,k+1\}, \qquad \widetilde r = \min\{r,k+2\}
\end{equation}
enter the error bounds. Since $H^{s}(\O) \subset H^{\widetilde s}(\O)$ with continuous embedding (and similarly for $\sigma$, $q$ and $r$), Assumption~\ref{ass:smoothness} also holds with $s$, $\sigma$, $q$, $r$ replaced by $\widetilde s$, $\widetilde \sigma$, $\widetilde q$, $\widetilde r$. Therefore, in the proofs below, we can, without loss of generality, assume that $s, \sigma, q \leq k+1$ and $r \leq k+2$.
 By the Sobolev embedding theorem, as $H^1(0,T;H^s(\O)) \subset C^0([0,T];H^s(\O))$ and $H^s(\O) \subset L^\infty(\O)$ for $s > 1$, Assumption \ref{ass:smoothness} immediately yields
\begin{equation}\label{eq:smooth-Linfty}
	u \in L^\infty(0,T;L^\infty(\O)), \qquad \o \in L^\infty(0,T;L^\infty(\O)), \qquad \fp \in L^\infty(0,T;L^\infty(\O)^2),
\end{equation}
the last inclusion holding since $r-1 > 1$ by \eqref{eq:smooth-indices}.
\begin{remark}
The indices $s$ and $\sigma$ are not independent: as $\fu = -\nabla u - \chit f(u) \fp$, if $f$ and $\fp$ are sufficiently smooth one can take $\sigma = s - 1$ (see the remark following \eqref{4.14dt}).
\end{remark}

We can now state the main result of the paper. 
By Theorem~\ref{thm:existence} the discrete solution exists on $[0,\Tmax)$; all bounds below are first proved for $t \leq \hTh \leq \Tmax$ (see the definition of $\hTh$ after Lemma~\ref{lemma4.1}) and extended to $[0,T]$ by Lemma~\ref{lem:continuation} and Corollary~\ref{cor:global}, consistently with  Remark \ref{rem:4.10}. We have the following theorem.
\begin{theorem}\label{thm:error}Let Assumption \ref{ass:smoothness} hold, let $\tu \simeq \tp \simeq 1$, and let the discrete initial data $u_h(0)$, $\p_h(0)$, $\o_h(0)$ be the $L^2(\O)$ orthogonal projections of $u_0$, $\p_0$, $\o_0$ onto $\Poly{k}(\Th)$. Then there exists $h_0$ such that for all $h<h_0$ the following error bounds hold:
	\begin{gather}
	\sup_{t\leq T}	{\| u(t) - u_h(t) \|^2_\O + \| \o(t) -  \o_h(t) \|_{\O}^{2}}\lesssim C(u)h^{2\widetilde s} + C(\fu)h^{2\widetilde\sigma} + C(\p)h^{2(\widetilde r-1)} + C(\o)h^{2\widetilde q} \label{eq:thm1}\\
		\int_0^T \left(\| \fu - \fu_h \|_{\O}^2	+ h \| u_h - \hu_h \|_{\Eh}^2\right) \lesssim
	C(u)h^{2\widetilde s} + C(\fu)h^{2\widetilde\sigma} + C(\p)h^{2(\widetilde r-1)} + C(\o)h^{2\widetilde q}  \label{eq:thm2}
	\end{gather}
	\begin{multline}\label{boundPhi}
		\| \p - \p_h \|_{L^\infty(0,T;L^2(\O))}^2
		+  \| \fp - \fp_h \|_{L^2(0,T;L^2(\O))}^2  +  \tp \| \p_h - \hp_h \|_{L^2(0,T;L^2(\Eh))}^2\\
		\lesssim
	C(u)h^{2\widetilde s} + C(\fu)h^{2\widetilde\sigma} + C(\p)h^{2(\widetilde r-1)} + C(\o)h^{2\widetilde q}
	\end{multline}
	with
	\begin{gather*}
	C(u) = \| u \|^2_{H^1(0,T;H^s(\O))}, \quad C(\fu) = \| \fu \|^2_{H^1(0,T;H^\sigma(\O))}, \\
	 C(\p)=\| \p \|^2_{H^1(0,T;H^r(\O))}, \quad C(\o) = \| \o \|^2_{L^\infty(0,T;{H^q(\O)})},
	\end{gather*}
	and with the capped indices $\widetilde s$, $\widetilde \sigma$, $\widetilde r$, $\widetilde q$ defined in \eqref{eq:capped-indices}.
\end{theorem}

\begin{corollary}
		 If $u \in H^1(0,T;H^{k+2}(\O))$, $\fu \in H^1(0,T;H^{k+1}(\O)^2)$, $\p \in H^1(0,T;H^{k+2}(\O))$ and $\o \in H^1(0,T;H^{k+1}(\O))$ all the variables converge with the optimal order $k+1$.
\end{corollary}

\begin{proof}
In order to prove an error bound, we rewrite equations \eqref{MS-HDG1}--\eqref{MS-HDG5} in a compact form. We introduce bilinear forms $\Au: \V \times \V \to \mathbb{R}$ and 
$\Ap: \V \times \V \to \mathbb{R}$, where  $\V = H^1(\Th) \times L^2(\Eh) \times (H^1(\Th))^2$, defined as follows: for
$U = (u,\hu,\fu) \in \V$, $W = (w,\hw,\fw) \in \V$
\begin{multline*}
\Au(U,W) = (\hfu \cdot \n, w)_{\Eh} - 
(\fu,\nabla w)_{\Th} + (\fu ,\fw)_{\Omega} 	- (u ,
\nabla \cdot \fw)_{\Th} \\+  (\hu , \fw\cdot \n)_{\Eh} - ( \underbrace{\fu\cdot\n + \tu (u - \hu)}_
{\hfu  \cdot \n}, \hw)_{\Eh},
\end{multline*}
  and for $\Phi = (\p,\hp,\fp) \in \V$, $\Theta = (\t,\ct,\ft) \in \V$ 
\begin{multline*}
\Ap(\Phi,\Theta) = ( \hfp\cdot\n, \t)_{\Eh}  -	(  \fp,\nabla\t)_{\Th}   + (\p,\t)_\Omega +  (\inu \fp,\ft)_{\Omega}
\\- (\p ,
\nabla \cdot \ft)_{\Th} + (\hp , \ft\cdot \n)_{\Eh} -  (\underbrace{\fp\cdot\n + \tp(\p - \hp)}_{\hfp \cdot \n},\ct)_{\Eh}.
\end{multline*}
We have the following result \cite{BBPP}.
\begin{proposition}
	It holds that
	\begin{gather}
	\Au(U,U) = \| \fu \|_{\Omega}^2 + \tu \|u - \hu \|_{\Eh}^2=: | U |_{\Au}^2,\\
		\Ap(\Phi,\Phi) =  \inu \| \fp \|_{\Omega}^2 +  \tp \| \p - \hp \|_{\Eh}^2  +  \| \p \|^2_{\O} =: \| \Phi \|_{\Ap}^2.
	\end{gather}
\end{proposition}
Both identities follow from an element by element integration by parts and from the definition of the numerical fluxes (see \cite{BBPP}, and the proof of Theorem \ref{thm:existence} in Section~\ref{sec:existence}, where the same computation is carried out for the discrete solution).
We also introduce the following compact notation for the nonlinear terms:
\begin{gather*}
\C \fp u \fw=  (\chit f(u) \fp ,\fw)_{\Omega}, \qquad
\B(u,w) = (u(u^{a-1}-1),w)_\O,\\
\G u \o \v = (g(u)(\o-1),\v) 
\end{gather*}
With this notation, letting $\Vh = \Sh = \Poly k(\Th) \times \Poly k(\Eh) \times (\Poly{k}(\Th))^2$ (we keep two different symbols for the spaces of the cellular and of the chemical unknowns, as they play different roles and as this allows a direct comparison with \cite{BBPP}, where the two spaces differ in the \HDGA\ case), the discrete problem rewrites as: find $U_h = (u_h,\hu_h,\fu_h) \in \Vh$, $\o_h \in \Poly{k}(\Th)$ and $\Phi_h = (\p_h,\hp_h,\fp_h)\in \Sh$ such that for all $W = (w,\hw,\fw) \in \Vh$, $\v \in \Poly{k}(\Th)$ and $\Theta = (\t,\ct,\ft)\in \Sh$ it holds that
\begin{gather}\label{MScompact1}
(\dt u_h,w) + (\dt \o_h, \v) + \Au(U_h,W) + \C{\fp_h} {u_h} \fw + \B(u_h,w) + \G {u_h}{\o_h} \v= (f_u,w)_\O\\ \label{MScompact2}
(\dt \p_h,\theta) + \Ap(\Phi_h,\Theta) - (\beta u_h + \delta \o_h,\theta)_\O = (f_\p,\theta)_\O.
\end{gather}
We start by observing that equation \eqref{MScompact2} is almost identical to the corresponding equation in the Keller-Segel HDG discretization method studied in \cite{BBPP}, from which it differs only in the term that couples it with equation \eqref{MScompact1}, in that  the bilinear form $(\beta u,\t)_\O$, is replaced here by $(\beta u + \delta \o,\t)_\O$. The following Lemma is therefore a straightforward modification of \cite[Lemma 4.8]{BBPP}.
\begin{lemma}\label{lemma4.1} If $\phi \in H^1(0,T;H^r(\Th))$, $2 < r \le k+2$, and $\p_h(0)$ is equal to the $L^2$ projection of $\p_0$, for all $T'\leq \Tmax$ it holds that
	\begin{multline}\label{4.5}
		\| \p - \p_h \|_{L^\infty(0,T';L^2(\O))}^2 
		+  \| \fp - \fp_h \|_{L^2(0,T';L^2(\O))}^2  +  \tp \| \p_h - \hp_h \|_{L^2(0,T';L^2(\Eh))}^2\\
		\lesssim	 
		h^{2(r-1)} \| \p \|^2_{H^1(0,T';H^{r}(\Th))}
		+	\| u- u_h \|_{L^2(0,T';L^2(\O))}^2 + \| \o - \o_h \|_{L^2(0,T';L^2(\O))}^2
	\end{multline}
\end{lemma}

Let $\Tmax \in (0,T]$ denote the maximal existence time of the discrete solution provided by Theorem \ref{thm:existence}, and set
\[
\upper := \sup_{t \leq T} \left( \| \o(t) \|_{L^\infty(\O)} + \| u(t) \|_{L^\infty(\O)} \right) + 2,
\]
which is finite under the regularity assumptions of Theorem \ref{thm:error}. We now let
\[
\hTh := \sup \Big\{ t' \in [0,\Tmax) :\ \| \o_h(t) \|_{L^\infty(\O)} + \| u_h(t) \|_{L^\infty(\O)} \leq \upper \ \text{ for all } t \leq t' \Big\}
\]
denote the largest time up to which the discrete solution satisfies the $L^\infty$ bound $\upper$ (with the convention $\hTh = 0$ if the set is empty, and $\hTh = T$ if $\Tmax = T$ and the bound holds on $[0,T]$). Reasoning as in \cite[Section 4.5]{BBPP}, we observe that $\|u_0 - u_h(0)\|_{L^\infty(\O)} \lesssim h^{s-1}\|u_0\|_{H^s(\O)}$ and $\|\o_0 - \o_h(0)\|_{L^\infty(\O)} \lesssim h^{q-1}\|\o_0\|_{H^q(\O)}$ ($s,q>1$) for the $L^2$ projected initial data (see, e.g., \cite{ern2021finite}) implies that for $h$ sufficiently small we have that $\hTh > 0$. We also observe that, if $\hTh < T$, then the discrete solution is defined on $[0,\hTh]$ and beyond: indeed, on $[0,\hTh)$ the function $u_h$ is bounded in $L^\infty(\O)$, hence in $L^2(\O)$, so that, by Theorem \ref{thm:existence} (i), either $\Tmax = T > \hTh$, or $\Tmax < T$ and $\| u_h(t) \|_\O$ is unbounded as $t \to \Tmax$, which forces $\hTh < \Tmax$. In both cases $\hTh < \Tmax$, and, by continuity in time of $u_h$ and $\o_h$, we have
\begin{equation}\label{eq:hThattained}
	\| \o_h(\hTh) \|_{L^\infty(\O)} + \| u_h(\hTh) \|_{L^\infty(\O)} = \upper \qquad \text{ if } 0 < \hTh < T.
\end{equation}
All the bounds below are stated for $t \leq \hTh$, where the discrete solution exists, and where $\| u_h(t) \|_{L^\infty(\O)} + \| \o_h(t) \|_{L^\infty(\O)} \leq \upper$.

We have the following Lemma, whose proof we postpone.
\begin{lemma}\label{lem:boundd} For all $t \leq \hTh$ it holds that
\begin{equation}\label{boundo}
	\| \o(t) - \o_h(t) \|_{\O}^2 \lesssim 
	h^{2q} \| \o \|_{L^\infty(0,T;H^q(\O))}^2 +  
	\int_0^t \| u - u_h \|_\O^2+  \int_0^t \| \o - \o_h  \|_\O^2,
\end{equation}
the hidden constant depending on the data but not on $h$.
\end{lemma}

\newcommand{\Cnorm}[2]{\| #2, \dt #2 \|_{#1} }
We now let $U^I = (u^I,\hu^I,\fu^I) \in \Vh$ be defined by
\begin{equation}\label{defUI}
(u^I,w)_\O + \Au(U^I,W)  = (u,w)_\O + \Au(U,W).
\end{equation}

 We know (see \cite[Lemma 4.6, case $k'=k$]{BBPP}, applied with $\fp = 0$ and $M=1$) that  the following bounds hold: if $u \in H^1(\O) \cap H^s(\Th)$, $\fu \in H^\sigma(\Th)^2$, $1 \leq s \leq k+1$, $1 \leq \sigma \leq k+1$ and $\hu = u|_{\Eh}$ then
	\begin{gather}\label{A.30}
		\|  \fu - \fuI \| _\O +
		h^{1/2} \| \uI - \huI \|_{\Eh} \lesssim h^{s} \| u \|_{H^s(\O)} + h^{\sigma} \| \fu \|_{H^\sigma(\O)},\\
\| u - \uI \|_{\O} \lesssim
			h^{s } \| u \|_{H^{s}(\O)} + h^{\sigma} \| \fu \|_{H^\sigma(\O)}, \label{4.12}\\
\label{4.14dt}		\| \dt u - \dt \uI \|_{\O}  \lesssim
			h^{s } 	\big(\| u \|_{H^s(\O)} + \| \dt u \|_{H^s(\O)}\big)  + h^{\sigma} \big(\| \fu \|_{H^\sigma(\O)} + \| \dt \fu \|_{H^\sigma(\O)}\big).
			\end{gather}
 We recall that, in \cite{BBPP}, these bounds are obtained by comparing $U^I$ with the stabilization-parametric HDG projection of $U$ introduced in \cite{cockburn2010projection}, which is possible since all the components of $\Vh$ are of the same degree $k$ and since $\tu$ does not depend on $h$; this is the only point where the analysis of the equal-order method differs from the one of the \HDGA\ method, for which the bound on $u - \uI$ requires a duality argument.

 Remark that, as $\fu = -\nabla u - \chit f(u) \fp$, if $u \in H^1(\O) \cap H^s(\Th)$ and $f(u) \fp \in H^{s-1}(\Th)^2$ (which holds provided $f$ and $\fp$ are sufficiently smooth), we can take $\sigma = s-1$, and the right hand sides of \eqref{A.30}--\eqref{4.14dt} are then of order $h^{s-1}$. Differently from the \HDGA\ case, where the two terms in the right hand side balance, in the equal-order case the convergence rate of all the variables is dictated by the regularity of the cellular flux: in order to obtain the optimal rate $h^{k+1}$ one needs $\fu \in H^{k+1}(\Th)^2$, which, for the above expression of $\fu$, holds provided $u \in H^{k+2}(\Th)$ and $f(u) \fp \in H^{k+1}(\Th)^2$. In the following we let
 \[
 \rho_h := h^{s} \big(\| u \|_{H^s(\O)} + \| \dt u \|_{H^s(\O)}\big) + h^{\sigma} \big(\| \fu \|_{H^\sigma(\O)} + \| \dt \fu \|_{H^\sigma(\O)}\big)
 \]
 denote the right hand side of \eqref{4.14dt}, and we observe that, by \eqref{4.12}, we also have $\| u - \uI \|_\O \lesssim \rho_h$.

 \renewcommand{\Zeta}{E}
 \renewcommand{\fzeta}{\feta}
\renewcommand{\b}{******}

We can now bound $U_h-U^I$. 
Let  $\xi = u_h - \uI$, $\hzeta = \hu_h - \huI$, $\fzeta = \fu_h - \fuI$, and let $\Zeta = (\xi,\hzeta,\fzeta)$.
We have
\begin{multline}
\frac 1 2 \dt \| \xi \|_\O^2 + | E |_{\Au}^2 = 
(\dt (u_h - \uI),\xi) +  \Au(U_h - U^I,E) \\
=
(\dt u_h,\xi) + \Au(U_h,E) - (\dt \uI,\xi) - \Au(U^I,E) 
\pm (\dt u,\xi) 
\pm \Au(U,E)
\\	\pm \C {\fp_h}{u_h}\feta
\pm \C {\fp}{u_h}\feta
\pm \C {\fp}{\uI}\feta  \pm \C {\fp}{u}\feta
\\	 \pm \B(u_h,\xi) \pm \B(\uI,\xi) \pm \B(u,\xi)
\\= \underbrace{	(\dt u_h,\xi) + \Au(U_h,E) + \C {\fp_h}{u_h}\feta + \B(u_h,\xi)}_{=(f_u,\xi)_\O}\\
\underbrace{- (\dt u,\xi) -   \Au(U,E)-  \C {\fp}{u}\feta -  \B(u,\xi) }_{=-(f_u,\xi)_\O}
+\underbrace{(\dt u,\xi) - (\dt \uI,\xi)}_{\lesssim \rho_h \| \xi \|_{\O}}+
\underbrace{\Au(U,E) - \Au(U^I,E)}_{= (\uI-u,\xi)_\O}
\\ +\underbrace{ \C {\fp - \fp_h}{u_h}\feta	}_{\lesssim \chit \| f(u_h) \|_{\infty} \| \fp - \fp_h \|_\O \| \feta\|_\O}
+ \underbrace{\C {\fp}{u^I}\feta  -  \C {\fp}{u_h}\feta}_{\leq \chit \| \fp \|_\infty \| f(\uI) - f(u_h) \|_\O \| \feta \|_\O}
+ \underbrace{\C {\fp}{u}\feta- \C {\fp}{u^I}\feta }_{\leq \chit \| \fp \|_\infty \| f(\uI) - f(u) \|_\O \| \feta \|_\O}
\\	+ \underbrace{\B(\uI,\xi) - \B(u_h,\xi)}_{\leq \| b(\uI) - b(u_h) \|_\O \| \xi \|_\O}
+ \underbrace{\B(u,\xi)  -  \B(\uI,\xi) }_{\leq \| b(\uI) - b(u) \|_\O \| \xi \|_\O},
\end{multline}
where $b(y) = y|y|^{a-1} - y$ is the function defined in \eqref{eq:extension}, and where we used that for $t  \leq \hTh$, $\| f(u_h(t)) \|_\infty \lesssim 1$.
We have
\begin{gather}
\label{4.14}\| f(v) - f(w) \|^2_\O \lesssim  (\| v \|_{L^\infty}^{2\ell} + \| w \|_{L^\infty}^{2\ell}) \| v - w 
\|_{\O}^2,\\
\| b(v) - b(w) \|^2_\O \lesssim  (\| v \|_{L^\infty}^{2(a-1)} + \| w \|_{L^\infty}^{2(a-1)} + 2) \| v - w 
\|_{\O}^2,
\end{gather}
where we used Assumption~F to  obtain \eqref{4.14}.

Then, since for $t \leq \hTh$ and $h$ sufficiently small we have $\| u_h \|_{L^\infty} \lesssim 1$ and, by the inverse inequality \eqref{eq:inverse}, the $L^\infty$ stability of the $L^2(\O)$ projection onto $\Poly{k}(\Th)$ and \eqref{4.12}, $\| \uI \|_{L^\infty} \lesssim \| u \|_{L^\infty} + h^{s-1} \| u \|_{H^s(\O)} + h^{\sigma-1} \| \fu \|_{H^\sigma(\O)} \lesssim 1$, it holds that
\begin{multline}
\frac 1 2 \dt \| \xi \|_\O^2 + | E |_{\Au}^2 \lesssim 
{ \rho_h \| \xi \|_{\O}} +
\| \uI-u\|_\O \| \xi \|_{\O}
+  \| \fp - \fp_h \|_\O \| \feta\|_\O \\
+ { \| \xi \|_\O \| \feta \|_\O}
+  \| \uI - u \|_\O \| \feta \|_\O
	+ \| \xi\|_\O^2
+  \|\uI - u \|_\O \| \xi \|_\O\\
\lesssim
\rho_h^2 + \| \xi \|_{\O}^2 + \frac 1 \vareps\| \fp - \fp_h \|_{\O}^2 + \vareps \| \feta \|_{\O}^2,
\end{multline}
where we used Young's inequality. Choosing $\vareps$ small enough and using $\| \feta \|_\O^2 \leq | E |^2_{\Au}$, the last term can be absorbed in the left hand side, and we obtain
\begin{equation}
\frac 1 2 \dt \| \xi \|_\O^2 + | E |_{\Au}^2 \lesssim
\rho_h^2  + \| \fp - \fp_h \|_{\O}^2 + \| \xi \|_{\O}^2,
\end{equation}
where, here and in the rest of the proof, the constants hidden in $\lesssim$ depend on $\upper$, $\| f \|_{W^{1,\infty}(0,\upper)}$, $a$, $\| \fp \|_{L^\infty(0,T;L^\infty(\O))}$ and $\chit$, but not on $h$,
and, integrating in $[0,t]$ and using  Lemma \ref{lemma4.1}
\begin{multline}\label{4.21}
\| \xi(t) \|_\O^2 + \int_0^t | E |_{\Au}^2	 \lesssim  \| \xi(0) \|_\O^2 + 
\int_0^t  \rho_h^2  + \int_0^t \| \fp - \fp_h \|_{\O}^2 + \int_0^t \| \xi \|_{\O}^2 
\\	\lesssim \| \xi(0) \|_\O^2 + 
h^{2s} \| u \|^2_{H^1(0,T;H^s(\O))} + h^{2\sigma} \| \fu \|^2_{H^1(0,T;H^\sigma(\O))} + h^{2(r-1)} \| \phi \|^2_{H^1(0,T;H^r(\O))} \\ + \int_0^t \| u - u_h \|_{\O}^2 + \int_{0}^t \| \o - \o_h \|^2_{\O},
\end{multline}
where the term $\int_0^t \| \xi \|^2_\O$ is absorbed through $\| \xi \|_\O \leq \| u - u_h \|_\O + \| u - \uI \|_\O \lesssim \| u - u_h \|_\O + \rho_h$. Since $u_h(0)$ is the $L^2(\O)$ orthogonal projection of $u_0$ onto $\Poly{k}(\Th)$, by \eqref{4.12} we have $\| \xi(0) \|_\O \leq \| u_h(0) - u_0 \|_\O + \| u_0 - \uI(0) \|_\O \lesssim h^{s} \| u_0 \|_{H^s(\O)} + h^{\sigma} \| \fu(0) \|_{H^\sigma(\O)}$, and a triangle inequality yields the following bound for $u - u_h$:
\begin{multline*}
\| u(t) - u_h(t) \|^2_\O \lesssim h^{2s} (\| u_0 \|^2_{H^s(\O)} + \| u \|^2_{H^1(0,T;H^s(\O))}) + h^{2\sigma} (\| \fu(0) \|^2_{H^\sigma(\O)} + \| \fu \|^2_{H^1(0,T;H^\sigma(\O))}) \\ + h^{2(r-1)}  \| \phi \|^2_{H^1(0,T;H^r(\Th))} + \int_0^t \| u - u_h \|_{\O}^2 + \int_{0}^t \| \o - \o_h \|^2_{\O}.
\end{multline*}
Adding \eqref{boundo} then yields
\begin{multline*}
\| u(t) - u_h(t) \|^2_\O + \| \o(t) -  \o_h(t) \|_{\O}^{2}\lesssim h^{2s}C(u) + h^{2\sigma}C(\fu) + h^{2(r-1)}C(\phi) + h^{2q}C(\o) \\ + \int_0^t (\| u - u_h \|_{\O}^2 + \| \o - \o_h \|^2_{\O}).
\end{multline*}

Gronwall's lemma finally yields, for all $t \leq \hTh$,
\begin{equation}\label{eq:errorbound}
{\| u(t) - u_h(t) \|^2_\O + \| \o(t) -  \o_h(t) \|_{\O}^{2}}\lesssim h^{2s}C(u) + h^{2\sigma}C(\fu) + h^{2(r-1)}C(\phi) + h^{2q}C(\o),
\end{equation}
with an implicit constant depending on  $\upper$ and $T$ but not on $h$.
Combining \eqref{eq:errorbound} with \eqref{4.21}, and using \eqref{A.30} and a triangle inequality (observe that, by \eqref{A.30}, $h \| \uI - \huI \|^2_\Eh \lesssim \rho_h^2$, while, by \eqref{4.21} and $\tu \simeq 1$, $\int_0^{\hTh} ( \| \fzeta \|^2_\O + \| \xi - \hzeta \|^2_\Eh) \lesssim \int_0^{\hTh} | E |^2_{\Au}$) yields 
\[
	\int_0^{\hTh} \left(\| \fu - \fu_h \|_{\O}^2	+ h \| u_h - \hu_h \|_{\Eh}^2\right) \lesssim
C(u)h^{2\widetilde s} + C(\fu)h^{2\widetilde\sigma} + C(\p)h^{2(\widetilde r-1)} + C(\o)h^{2\widetilde q} 
\] and Lemma \ref{lemma4.1} yields 
\begin{multline*}
	\| \p - \p_h \|_{L^\infty(0,\hTh;L^2(\O))}^2 
+  \| \fp - \fp_h \|_{L^2(0,\hTh;L^2(\O))}^2  +  \tp \| \p_h - \hp_h \|_{L^2(0,\hTh;L^2(\Eh))}^2\\
\lesssim	 
C(u)h^{2\widetilde s} + C(\fu)h^{2\widetilde\sigma} + C(\p)h^{2(\widetilde r-1)} + C(\o)h^{2\widetilde q} .
\end{multline*}
We conclude the proof thanks to the following Lemma.
\begin{lemma}\label{lem:continuation}
	There exists $h^*$ such that for $h\leq h^*$ we have $\hTh = T$.
\end{lemma}
Indeed, the bounds above hold in $[0,\hTh]$ with constants independent of $h$. For $h \leq h^*$ this gives  \eqref{eq:thm1}--\eqref{boundPhi}.
\end{proof}

\newcommand{\Tstarhn}{\widehat T_{h_n}}
\newcommand{\psmooth}{q}

\begin{corollary}[Global existence of the discrete solution for $h$ small]\label{cor:global}
	Under the assumptions of Theorem \ref{thm:error}, let $h^*$ be given by Lemma \ref{lem:continuation}. Then, for all $h \leq h^*$, the unique discrete solution of Problem \eqref{MS-HDG1}--\eqref{MS-HDG5} given by Theorem \ref{thm:existence} is defined on the whole interval $[0,T]$, that is $\Tmax = T$, and it satisfies
	\[
	\sup_{t \leq T} \left( \| \o_h(t) \|_{L^\infty(\O)} + \| u_h(t) \|_{L^\infty(\O)} \right) \leq \upper.
	\]
\end{corollary}
\begin{proof}
	By Lemma \ref{lem:continuation}, for $h \leq h^*$ we have $\hTh = T$, so that $\| u_h(t) \|_{L^\infty(\O)} + \| \o_h(t) \|_{L^\infty(\O)} \leq \upper$ for all $t < \Tmax$. In particular $\sup_{t < \Tmax} \| u_h(t) \|_\O \leq |\O|^{1/2} \upper < \infty$, and Theorem \ref{thm:existence} (i) yields $\Tmax = T$ and the extension of the solution to $[0,T]$; the $L^\infty$ bound at $t = T$ follows by continuity.
\end{proof}

\begin{remark}\label{rem:4.10}
	Theorem \ref{thm:existence} and Corollary \ref{cor:global} justify the strategy outlined in the paragraph preceding Theorem \ref{thm:error}, which ensures that, for $h$ small enough, the discrete solution exists on the whole interval $[0,T]$: the error bounds are first established on $[0,\hTh]$, where the solution is known to exist and to be bounded in $L^\infty(\O)$, and Lemma \ref{lem:continuation} then shows that $\hTh = T$ for $h \leq h^*$. Notice that no positivity of the discrete solution is required at any stage.
\end{remark}

We now prove Lemmas \ref{lem:boundd} and \ref{lem:continuation}.
\begin{proof}[Proof of Lemma \ref{lem:boundd}] Let $\o^\pi \in \Poly{k}(\Th)$ denote the $L^2$ orthogonal projection of $\o$. We have that $\o \in H^{q}(\O)$, $0 \leq q \leq k+1$ implies
	\[
	\| \o - \o^\pi \|_{\O} \lesssim h^{q} | \o |_{q,\O}.
	\]
	Testing \eqref{MS2}, that is $\dt \o + g(u)\o = g(u)$, with $\v \in \Poly{k}(\Th)$, we can write
	\begin{multline*}
		\dt(\o^\pi,\v)_\O + (g(u)\o^\pi,\v)_\O = \dt(\o,\v)_\O + (g(u)\o,\v)_\O + \dt(\o^\pi - \o,\v)_\O + (g(u)(\o^\pi - \o),\v)_\O \\ =  (g(u),\v)_\O + \dt(\o^\pi - \o,\v)_\O + (g(u)(\o^\pi - \o),\v)_\O,
	\end{multline*}
	where $\dt(\o^\pi - \o,\v)_\O = 0$, since $\dt \o^\pi = (\dt \o)^\pi$ and $\v \in \Poly{k}(\Th)$.
	Let now $\mathfrak e = \o^\pi - \o_h$. We have
	\begin{multline*}
	(	\dt{\mathfrak e},{\mathfrak e})_\O + (g(u)\o^\pi,{\mathfrak e}) \pm(g(u_h) \o^\pi,{\mathfrak e})_\O - (g(u_h)\o_h,{\mathfrak e})_\O \\ = 
	(	\dt{\mathfrak e},{\mathfrak e})_\O + (g(u_h){\mathfrak e},{\mathfrak e})_\O - ((g(u_h) - g(u))\o^\pi,{\mathfrak e})_\O 
		\\= (g(u)-g(u_h),{\mathfrak e})_\O +  \dt(\o^\pi - \o,{\mathfrak e})_\O + (g(u)(\o^\pi - \o),{\mathfrak e})_\O,
	\end{multline*}
	that is
	\begin{multline}\label{4.5b} \frac12 \dt \| {\mathfrak e} \|_\O^2 \leq
	\frac12	\dt \| {\mathfrak e} \|_\O^2 + (g(u_h){\mathfrak e},{\mathfrak e}) \\ = ((g(u_h) - g(u) )(\o^\pi - 1),{\mathfrak e})_\O +  \dt\underbrace{(\o^\pi - \o,\e)_\O }_{=0\text{ as }\e\in \Poly{k}}+ (g(u)(\o^\pi - \o),{\mathfrak e})_\O\\
		\lesssim (\| \o^\pi \|_{L^\infty(\O)} + 1) \| g(u_h) - g(u) \|_\O \| {\mathfrak e} \|_{\O} + \| g \|_{L^\infty} \| \o - \o^\pi \|_{\O} \| \e \|_\O \\ \lesssim
		(\| \o^\pi \|_{L^\infty(\O)} + 1)  \| g \|_{W^{1,\infty}(0,\upper)} \| u_h - u \|_\O \| {\mathfrak e} \|_{\O} + \| g \|_{L^\infty(0,\upper)} h^{q} | \o |_{q,\O} \| \e \|_\O ,
	\end{multline}
	where we used that, $\Th$ being quasi-uniform, the $L^2(\O)$ orthogonal projection onto $\Poly{k}(\Th)$ is stable in $L^\infty(\O)$, so that $\| \o^\pi \|_{L^\infty(\O)} \lesssim \| \o \|_{L^\infty(\O)} \lesssim 1$ by \eqref{eq:smooth-Linfty}. Moreover, for $t \leq \hTh$ both $u$ and $u_h$ take values in $[-\upper,\upper]$, and, thanks to the extension \eqref{eq:extension}, $g$ is Lipschitz continuous on $[-\upper,\upper]$ with constant $\max_{[0,\upper]} | g' | \leq \| g \|_{W^{1,\infty}(0,\upper)}$, which gives $\| g(u_h) - g(u) \|_\O \leq \| g \|_{W^{1,\infty}(0,\upper)} \| u_h - u \|_\O$.
	Applying Young inequality and integrating in $[0,t]$ yields
	\begin{equation*}\| {\mathfrak e}(t) \|^2_\O 
		\lesssim  \| {\mathfrak e}(0) \|^2_\O +  h^{2q} \| \o \|_{L^2(0,T;H^q(\O))}^2  + 
		\int_0^t \| u - u_h \|_\O^2 +  \int_0^t \| {\mathfrak e} \|_\O^2.
	\end{equation*}
	Since $\o_h(0)$ is the $L^2(\O)$ orthogonal projection of $\o_0$, we have $\e(0) = 0$. Moreover, $\| \o \|_{L^2(0,T;H^q(\O))} \leq T^{1/2} \| \o \|_{L^\infty(0,T;H^q(\O))}$, $\int_0^t \| \e \|_\O^2 \lesssim \int_0^t \| \o - \o_h \|_\O^2 + h^{2q} T \| \o \|^2_{L^\infty(0,T;H^q(\O))}$ and $\| \o(t) - \o^\pi(t) \|_\O \lesssim h^q \| \o \|_{L^\infty(0,T;H^q(\O))}$. A triangle inequality then yields \eqref{boundo}.
\end{proof}

\begin{proof}[Proof of Lemma \ref{lem:continuation}]
We prove the result by contradiction.
If the statement of the lemma is false, there exists a sequence $h_n \to 0$ such that $\Tstarhn < T$.	
By definition of $\Tstarhn$, since $\Tstarhn < T$, the discrete solution exists on a neighbourhood of $\Tstarhn$ (see the discussion preceding Lemma \ref{lem:boundd}). Moreover, by the same discussion, $\Tstarhn > 0$ for $n$ large enough, so that \eqref{eq:hThattained} holds with $h = h_n$, that is
\begin{equation}\label{eq:contradicthp}
	\| u_{h_n}(\Tstarhn) \|_{L^\infty(\O)} + \| \o_{h_n}(\Tstarhn) \|_{L^\infty(\O)} = \upper.
\end{equation}
Thanks to \eqref{eq:errorbound} we also have
\begin{multline*}
	\|  u_{h_n}(\Tstarhn) \|_{L^\infty(\Omega)} + 
	\| \o_{h_n}(\Tstarhn) \|_{L^\infty(\Omega)} 
	  \leq  \| u(\Tstarhn) \|_{L^\infty(\O)} +  \| \o(\Tstarhn) \|_{L^\infty(\O)} \\
+	   \|  u_{h_n}(\Tstarhn) - u(\Tstarhn)  \|_{L^\infty(\Omega)} 
	  + \|  \o_{h_n}(\Tstarhn) - \o(\Tstarhn)  \|_{L^\infty(\Omega)}  \\
  \leq   \| u(\Tstarhn) \|_{L^\infty(\O)}  +  \| \o(\Tstarhn) \|_{L^\infty(\O)} + C h_n^{\min\{s,\sigma,r-1,q\} - 1} ,
\end{multline*}
with $C$ independent of $n$, where, in order to bound the $L^\infty(\O)$ norms of $u_{h_n} - u$ and $\o_{h_n} - \o$, we added and subtracted the $L^2(\O)$ orthogonal projections of $u$ and $\o$ onto $\Poly{k}(\Th)$, and we used the inverse inequality \eqref{eq:inverse}, the bound \eqref{eq:errorbound}, and the approximation properties of the $L^2(\O)$ orthogonal projection, both in $L^2(\O)$ and in $L^\infty(\O)$ (see, e.g., \cite{ern2021finite}). As, under our smoothness assumptions, $\min\{s,\sigma,r-1,q\} > 1$, the right hand side is bounded by $\upper - 2 + C h_n^{\min\{s,\sigma,r-1,q\} - 1}$, which, for $n$ large enough, is strictly smaller than $\upper$. This contradicts \eqref{eq:contradicthp}.
\end{proof}

\section{Proof of Theorem \ref{thm:existence}}

\label{sec:existence}

Problem \eqref{MS-HDG1}--\eqref{MS-HDG5} is a system of ordinary differential equations for $u_h$, $\o_h$, $\p_h$, coupled with the algebraic equations \eqref{MS-HDG3}, \eqref{MS-HDG6}, \eqref{MS-HDG4}, \eqref{MS-HDG5} for the fluxes and the traces. In this section we show that it admits a unique solution, and we give conditions under which the solution is defined on the whole interval $[0,T]$. We stress that the argument is purely finite dimensional and does not rely on the positivity of $u_h$, which, differently from the continuous case (see Theorem \ref{existunique}), is not guaranteed by the discretization. We observe that, since $f \in C^2([0,\infty))$ with $f(0) = 0$, $g \in C^1([0,\infty))$ and $a > 1$, the extended $f$, $g$ and $b$ defined by \eqref{eq:extension} are locally Lipschitz continuous on $\mathbb R$, and they satisfy
\begin{equation}\label{eq:extprop}
	g(y) \geq 0, \qquad b(y)\, y = |y|^{a+1} - y^2 \geq - y^2 \qquad \text{ for all } y \in \mathbb R,
\end{equation}
and $|f(y)| \lesssim |y|^\expb$ for all $y \in \mathbb R$. We also let, for $\rho \geq 0$,
\[
L_f(\rho) = \max_{[0,\rho]} |f'|, \qquad L_g(\rho) = \max_{[0,\rho]} |g'|, \qquad L_b(\rho) = a \rho^{a-1} + 1,
\]
which are Lipschitz constants of $f$, $g$ and $b$ on $[-\rho,\rho]$.

In the following, the integrals defining the nonlinear terms $\C{\fp}{u}{\fw}$, $\B(u,w)$ and $\G{u}{\o}{\v}$ are assumed to be computed exactly. Moreover, we will repeatedly use the inverse inequality
\begin{equation}\label{eq:inverse}
	\| v \|_{L^\infty(\O)} \leq \Cinv h^{-1} \| v \|_\O \qquad \text{ for all } v \in \Poly{k}(\Th),
\end{equation}
which holds, with $\Cinv$ depending only on $k$ and on the shape regularity and quasi-uniformity constants of $\Th$ (recall that $\O \subset \mathbb R^2$).

We start by observing that the algebraic equations can be solved element-wise in terms of $u_h$ and $\p_h$.
\begin{lemma}[Local solvers]\label{lem:localsolvers}
	\begin{enumerate}[(i)]
		\item For every $\p_h \in \Poly{k}(\Th)$ there exists a unique pair $(\fp_h,\hp_h) \in (\Poly{k}(\Th))^2 \times \Poly{k}(\Eh)$ satisfying \eqref{MS-HDG4}--\eqref{MS-HDG5}. The map $\Lp: \p_h \mapsto (\fp_h,\hp_h)$ is linear.
		\item For every $u_h \in \Poly{k}(\Th)$ and every $\bs{r} \in (L^2(\O))^2$ there exists a unique pair $(\fu_h,\hu_h) \in (\Poly{k}(\Th))^2 \times \Poly{k}(\Eh)$ such that, for all $\fv \in (\Poly{k}(\Th))^2$ and $\hw \in \Poly{k}(\Eh)$,
		\begin{equation}\label{eq:localu}
			(\fu_h,\fv)_\O - (u_h,\nabla\cdot\fv)_\Th + (\hu_h,\fv\cdot\n)_\Eh = -(\bs r,\fv)_\O, \qquad
			(\fu_h\cdot\n + \tu(u_h - \hu_h),\hw)_\Eh = 0.
		\end{equation}
		The map $\Lu: (u_h,\bs r) \mapsto (\fu_h,\hu_h)$ is linear and continuous from $\Poly{k}(\Th) \times (L^2(\O))^2$ to $(\Poly{k}(\Th))^2 \times \Poly{k}(\Eh)$.
	\end{enumerate}
	In particular, $(u_h,\hu_h,\fu_h)$ and $(\p_h,\hp_h,\fp_h)$ satisfy \eqref{MS-HDG3}, \eqref{MS-HDG6}, \eqref{MS-HDG4}, \eqref{MS-HDG5} if and only if
	\begin{equation}\label{eq:algebraic}
		(\fp_h,\hp_h) = \Lp \p_h, \qquad (\fu_h,\hu_h) = \Lu(u_h, \chit f(u_h) \fp_h).
	\end{equation}
\end{lemma}
\begin{proof}
	In both cases the unknowns and the test functions range over spaces of the same (finite) dimension, so that we have a square linear system and it is sufficient to prove uniqueness, that is, to show that the homogeneous system only has the trivial solution. Let us consider (i) with $\p_h = 0$: taking $\ft = \fp_h$ in \eqref{MS-HDG4} and $\ct = \hp_h$ in \eqref{MS-HDG5}, and adding the two identities, we obtain
	\[
	\inu \| \fp_h \|_\O^2 + (\hp_h,\fp_h\cdot\n)_\Eh - (\fp_h\cdot\n,\hp_h)_\Eh + \tp \| \hp_h \|_\Eh^2 = \inu \| \fp_h \|_\O^2 + \tp \| \hp_h \|_\Eh^2 = 0,
	\]
	whence $\fp_h = 0$ and $\hp_h = 0$. The same argument, with $u_h = 0$ and $\bs r = 0$ in \eqref{eq:localu}, $\fv = \fu_h$ and $\hw = \hu_h$, gives $\| \fu_h \|_\O^2 + \tu \| \hu_h \|_\Eh^2 = 0$ and proves uniqueness in (ii). Linearity is obvious, and continuity of $\Lu$ follows from the fact that the right hand side of \eqref{eq:localu} depends continuously on $(u_h,\bs r)$ and that the solution operator of a fixed invertible linear system on a finite dimensional space is continuous. Finally, \eqref{eq:algebraic} is just a rewriting of \eqref{MS-HDG3}, \eqref{MS-HDG6}, \eqref{MS-HDG4}, \eqref{MS-HDG5}, since $\chit f(u_h)\fp_h \in (L^2(\O))^2$ for $u_h,\p_h \in \Poly{k}(\Th)$.
\end{proof}

\begin{proof}[Proof of Theorem \ref{thm:existence}]
	We split the proof in three steps.
	
	\medskip\noindent
	{\em Step 1: reduction to a system of ODEs and local existence.} Thanks to Lemma \ref{lem:localsolvers}, Problem \eqref{MS-HDG1}--\eqref{MS-HDG5} is equivalent to finding $(u_h,\o_h,\p_h) \in C^1([0,\Tmax);\Poly{k}(\Th)^3)$ such that, for all $w, \v, \t \in \Poly{k}(\Th)$ and for all $t$,
	\begin{gather}\label{eq:ode1}
		(\dt u_h,w)_\O = (f_u,w)_\O - (\hfu_h\cdot\n,w)_\Eh + (\fu_h,\nabla w)_\Th - \B(u_h,w),\\ \label{eq:ode2}
		(\dt \o_h,\v)_\O = - \G{u_h}{\o_h}{\v},\\ \label{eq:ode3}
		(\dt \p_h,\t)_\O = (f_\p,\t)_\O - (\hfp_h\cdot\n,\t)_\Eh + (\fp_h,\nabla\t)_\Th - (\p_h,\t)_\O + (\beta u_h + \delta \o_h,\t)_\O,
	\end{gather}
	where $(\fp_h,\hp_h)$ and $(\fu_h,\hu_h)$ are given by \eqref{eq:algebraic}, and $\hfu_h\cdot\n = \fu_h\cdot\n + \tu(u_h - \hu_h)$, $\hfp_h\cdot\n = \fp_h\cdot\n + \tp(\p_h-\hp_h)$. Since $(\cdot,\cdot)_\O$ is a scalar product on the finite dimensional space $\Poly{k}(\Th)$, equations \eqref{eq:ode1}--\eqref{eq:ode3} are equivalent to a system of ordinary differential equations
	\begin{equation}\label{eq:odesystem}
		\dt (u_h,\o_h,\p_h) = \mathsf{N}_h(t;u_h,\o_h,\p_h), \qquad (u_h,\o_h,\p_h)(0) = (u_h(0),\o_h(0),\p_h(0)),
	\end{equation}
	for a map $\mathsf N_h: [0,T] \times \Poly{k}(\Th)^3 \to \Poly{k}(\Th)^3$. We claim that $\mathsf N_h$ is continuous in $t$ and locally Lipschitz continuous with respect to $(u_h,\o_h,\p_h)$, uniformly in $t$. Continuity in $t$ only enters through the terms $(f_u,w)_\O$ and $(f_\p,\t)_\O$, and follows from $f_u, f_\p \in C^0([0,T];L^2(\O))$. All the terms in the right hand side of \eqref{eq:ode1}--\eqref{eq:ode3} are linear in $(u_h,\o_h,\p_h)$, $(\fp_h,\hp_h) = \Lp\p_h$ and $(\fu_h,\hu_h)$, except for $\B(u_h,w)$, $\G{u_h}{\o_h}{\v}$ and $\fu_h$, which depends on $u_h$ and $\p_h$ also through $\Lu(0,\chit f(u_h)\fp_h)$. Since linear maps between finite dimensional spaces are Lipschitz, and since compositions, sums and products of locally Lipschitz maps are locally Lipschitz, it is sufficient to check that the maps
	\[
	u_h \mapsto b(u_h) \in L^2(\O), \qquad (u_h,\o_h) \mapsto g(u_h)(\o_h - 1) \in L^2(\O), \qquad (u_h,\p_h) \mapsto f(u_h)\fp_h \in (L^2(\O))^2
	\]
	are locally Lipschitz on $\Poly{k}(\Th)$, resp. $\Poly{k}(\Th)^2$, endowed with the $L^2(\O)$ norm. Let $R > 0$ and let $v, w, \o_1, \o_2, \p_1, \p_2 \in \Poly{k}(\Th)$ with $L^2(\O)$ norm bounded by $R$. By \eqref{eq:inverse} their $L^\infty(\O)$ norms are bounded by $\rho_h := \Cinv h^{-1} R$, and, letting $\fp_i$ denote the first component of $\Lp \p_i$, $i = 1,2$, we have
	\begin{gather*}
		\| b(v) - b(w) \|_\O \leq L_b(\rho_h) \| v - w \|_\O, \\
		\| g(v)(\o_1 - 1) - g(w)(\o_2 - 1) \|_\O \leq L_g(\rho_h) \| v - w \|_{L^\infty(\O)} \| \o_1 - 1 \|_\O + \max_{[0,\rho_h]} |g|\, \| \o_1 - \o_2 \|_\O, \\
		\| f(v)\fp_1 - f(w)\fp_2 \|_\O \leq L_f(\rho_h) \| v - w \|_{L^\infty(\O)} \| \fp_1 \|_\O + \max_{[0,\rho_h]} |f|\, \| \fp_1 - \fp_2 \|_\O,
	\end{gather*}
	and the claim follows by using once again \eqref{eq:inverse} to bound $\| v - w \|_{L^\infty(\O)} \leq \Cinv h^{-1} \| v - w \|_\O$, and the boundedness of the linear map $\Lp$ to bound $\| \fp_1 \|_\O \lesssim R$ and $\| \fp_1 - \fp_2 \|_\O \lesssim \| \p_1 - \p_2 \|_\O$ (with constants depending on $h$). The Picard--Lindel\"of theorem then yields the existence of a unique maximal solution $(u_h,\o_h,\p_h) \in C^1([0,\Tmax);\Poly{k}(\Th)^3)$ of \eqref{eq:odesystem}, for some $\Tmax \in (0,T]$, with the following alternative: either $\Tmax = T$ and the solution extends to $[0,T]$, or
	\begin{equation}\label{eq:blowup}
		\limsup_{t \to \Tmax} \Big( \| u_h(t) \|_\O + \| \o_h(t) \|_\O + \| \p_h(t) \|_\O \Big) = + \infty.
	\end{equation}
	The regularity of $(\fp_h,\hp_h) = \Lp\p_h$ and of $(\fu_h,\hu_h) = \Lu(u_h,\chit f(u_h)\fp_h)$ stated in (i) follows from \eqref{eq:algebraic}, since $\Lp$ and $\Lu$ are linear and continuous and $t \mapsto f(u_h(t))$ is continuous with values in $L^\infty(\O)$. Uniqueness of the solution of Problem \eqref{MS-HDG1}--\eqref{MS-HDG5} follows from uniqueness for \eqref{eq:odesystem} and from Lemma \ref{lem:localsolvers}. To conclude the proof of (i) it remains to show that \eqref{eq:blowup} cannot occur if $\sup_{t<\Tmax} \| u_h(t) \|_\O < \infty$; this will follow from the bounds \eqref{eq:boundo_stab} and \eqref{eq:boundp_stab} that we prove in the next step.
	
	\medskip\noindent
	{\em Step 2: bounds on $\o_h$ and $\p_h$.} Since constants belong to $\Poly{k}(\Th)$, we can take $\v = \o_h - 1$ in \eqref{MS-HDG7}, obtaining
	\[
	\frac 12 \dt \| \o_h - 1 \|_\O^2 + (g(u_h)(\o_h - 1),\o_h - 1)_\O = 0,
	\]
	and, since $g \geq 0$ by \eqref{eq:extprop}, the second term is non negative, so that $\dt \| \o_h - 1 \|_\O^2 \leq 0$, which gives \eqref{eq:boundo_stab}. Let us now take $\t = \p_h$ in \eqref{MS-HDG2}, $\ft = \fp_h$ in \eqref{MS-HDG4} and $\ct = \hp_h$ in \eqref{MS-HDG5}, and add the resulting identities. Integrating by parts element by element we have $-(\fp_h,\nabla\p_h)_\Th - (\p_h,\nabla\cdot\fp_h)_\Th = -(\p_h,\fp_h\cdot\n)_\Eh$, and therefore
	\begin{multline*}
		\frac 12 \dt \| \p_h \|_\O^2 + (\hfp_h\cdot\n,\p_h - \hp_h)_\Eh - (\fp_h\cdot\n,\p_h - \hp_h)_\Eh + \| \p_h \|_\O^2 + \inu \| \fp_h \|_\O^2 \\
		= \frac 12 \dt \| \p_h \|_\O^2 + \| \p_h \|_\O^2 + \inu \| \fp_h \|_\O^2 + \tp \| \p_h - \hp_h \|_\Eh^2 = (\beta u_h + \delta \o_h + f_\p,\p_h)_\O,
	\end{multline*}
	where we used the definition of the numerical flux $\hfp_h\cdot\n = \fp_h\cdot\n + \tp(\p_h - \hp_h)$. By Young's inequality, $(\beta u_h + \delta \o_h + f_\p,\p_h)_\O \leq \| \p_h \|_\O^2 + \frac 34 (\beta^2 \| u_h \|_\O^2 + \delta^2 \| \o_h \|_\O^2 + \| f_\p \|_\O^2)$, whence
	\begin{equation}\label{eq:energyp}
		\dt \| \p_h \|_\O^2 + 2 \inu \| \fp_h \|_\O^2 + 2 \tp \| \p_h - \hp_h \|_\Eh^2 \leq \frac 32 \Big( \beta^2 \| u_h \|_\O^2 + \delta^2 \| \o_h \|_\O^2 + \| f_\p \|_\O^2 \Big),
	\end{equation}
	and \eqref{eq:boundp_stab} follows by integration in time. This proves (ii). Now, if $\sup_{t < \Tmax} \| u_h(t) \|_\O < \infty$, then \eqref{eq:boundo_stab} and \eqref{eq:boundp_stab} imply that also $\| \o_h(t) \|_\O$ and $\| \p_h(t) \|_\O$ are bounded on $[0,\Tmax)$, so that \eqref{eq:blowup} is excluded, $\Tmax = T$, and the solution extends to $[0,T]$. This concludes the proof of (i).
	
	\medskip\noindent
	{\em Step 3: bound on $u_h$ for bounded $f$.} We take $w = u_h$ in \eqref{MS-HDG1}, $\fv = \fu_h$ in \eqref{MS-HDG3} and $\hw = \hu_h$ in \eqref{MS-HDG6}, and add. Reasoning exactly as in Step 2 (element-wise integration by parts and definition of $\hfu_h\cdot\n$) we obtain
	\begin{equation}\label{eq:energyu}
		\frac 12 \dt \| u_h \|_\O^2 + \| \fu_h \|_\O^2 + \tu \| u_h - \hu_h \|_\Eh^2 + \B(u_h,u_h) + \C{\fp_h}{u_h}{\fu_h} = (f_u,u_h)_\O.
	\end{equation}
	By \eqref{eq:extprop}, $\B(u_h,u_h) = \int_\O b(u_h) u_h \geq - \| u_h \|_\O^2$. Letting $M_f := \| f \|_{L^\infty(0,\infty)} = \sup_{y \in \mathbb R} |f(y)| < \infty$, we have
	\[
	| \C{\fp_h}{u_h}{\fu_h} | \leq \chit M_f \| \fp_h \|_\O \| \fu_h \|_\O \leq \frac 12 \| \fu_h \|_\O^2 + \frac 12 \chit^2 M_f^2 \| \fp_h \|_\O^2,
	\]
	and $(f_u,u_h)_\O \leq \frac 12 \| f_u \|_\O^2 + \frac 12 \| u_h \|_\O^2$. Plugging these bounds in \eqref{eq:energyu} yields
	\begin{equation}\label{eq:energyu2}
		\dt \| u_h \|_\O^2 + \| \fu_h \|_\O^2 + 2 \tu \| u_h - \hu_h \|_\Eh^2 \leq 3 \| u_h \|_\O^2 + \| f_u \|_\O^2 + \chit^2 M_f^2 \| \fp_h \|_\O^2.
	\end{equation}
	We now set $\lambda := \chi^2 M_f^2$, so that $2 \lambda \inu \| \fp_h \|_\O^2 - \chit^2 M_f^2 \| \fp_h \|_\O^2 = \lambda \inu \| \fp_h \|_\O^2$, and we add \eqref{eq:energyu2} to $\lambda$ times \eqref{eq:energyp}. Letting $E(t) := \| u_h(t) \|_\O^2 + \lambda \| \p_h(t) \|_\O^2$ we get
	\begin{multline}\label{eq:energytot}
		\dt E + \| \fu_h \|_\O^2 + 2 \tu \| u_h - \hu_h \|_\Eh^2 + \lambda \inu \| \fp_h \|_\O^2 + 2 \lambda \tp \| \p_h - \hp_h \|_\Eh^2 \\
		\leq \Big( 3 + \frac 32 \lambda \beta^2 \Big) \| u_h \|_\O^2 + \| f_u \|_\O^2 + \frac 32 \lambda \Big( \delta^2 \| \o_h \|_\O^2 + \| f_\p \|_\O^2 \Big)
		\leq C_1 E + C_1 \Big( \| f_u \|_\O^2 + \| f_\p \|_\O^2 + \| \o_h \|_\O^2 \Big),
	\end{multline}
	with $C_1 = \max\{ 3 + \frac 32 \lambda \beta^2, 1, \frac 32 \lambda \delta^2, \frac 32 \lambda \}$. By \eqref{eq:boundo_stab}, $\| \o_h(t) \|_\O \leq |\O|^{1/2} + \| \o_h(0) - 1 \|_\O$ for all $t < \Tmax$. Gronwall's lemma then gives, for all $t < \Tmax$,
	\[
	E(t) \leq e^{C_1 t} \Big( E(0) + C_1 \int_0^t \big( \| f_u \|_\O^2 + \| f_\p \|_\O^2 \big) + C_1 t \big( |\O|^{1/2} + \| \o_h(0) - 1 \|_\O \big)^2 \Big) \leq C_2,
	\]
	with $C_2$ depending only on the quantities listed in the statement, and, in particular, not on $h$ and not on $t \in [0,T]$. In particular $\sup_{t<\Tmax} \| u_h(t) \|_\O < \infty$, so that, by (i), $\Tmax = T$. Integrating \eqref{eq:energytot} on $[0,T]$ and using again \eqref{eq:boundo_stab} finally yields \eqref{eq:stability}, which proves (iii).
\end{proof}

\newcommand{\code}[1]{\texttt{#1}}
\newcommand{\jj}{\bm{j}}
\newcommand{\psivec}{\bm{\psi}}
\newcommand{\nn}{\bm{n}}
\newcommand{\PP}{\mathbb{P}}
\newcommand{\dd}{\,\mathrm{d}}

\section{Numerical tests}\label{sec:numerical}

In this section, we present numerical experiments to test the performance of the proposed HDG method.
We consider two test cases: 
the first is a problem with known exact solution, used to verify the convergence rates of the method while the second shows how the method can reproduce the Turing patterns of the MS
model.

We test the equal-order formulation \HDGk\ introduced and analyzed in
Section~\ref{sec:HDG_MS}, for $k=1,2$, i.e.\ \HDGku\ and \HDGkd: the cellular
flux $\fu$, the density $u$, the chemoattractant $\p$, its flux $\fp$, the
traces $\hu$, $\hp$ and the density
$\o$ are all approximated by polynomials
of degree $k$.
Following \cite{lombardo2017demyelination} we take
\[
f(u)=\frac{u}{1+u},\qquad a=2,\qquad g(u)=\frac{u^2}{1+u}.
\]
All numerical experiments were performed in MATLAB R2026a 
on a MacBook Pro with an Apple  M5 Max processor and 128 GB of RAM.

\begin{remark}[Implementation]\label{rem:implementation}
	At each time step the local variables are eliminated by static condensation, so
	that the global linear system involves only the traces $\widehat u_h$ and
	$\widehat\varphi_h$. Since the equation for $\mathfrak d$ contains no spatial
	derivatives, $\mathfrak d_h$ is updated element by element and adds no globally
	coupled unknown. Time derivatives are discretized by the backward Euler method in
	Example~1 and by the BDF2 formula in Example~2, where long time intervals and
	large time steps are needed. In Example~1 the nonlinear terms are treated, as in
	\cite{BBPP}, by a Picard iteration in which the chemical flux $\boldsymbol\psi_h$
	is frozen at the previous iterate; with the small time step used there, it
	converges in three iterations per step. Its convergence however deteriorates as
	the time step grows: for $ dt=0.1$, as in Example~2, it requires about ten
	iterations per step. In that case we use instead a Newton linearization of the
	chemotactic term, which only modifies the local problems and converges in 
	about two or
	three iterations. 
	The iterations are stopped when the $L^2(\Omega)$ norm of the
	increment of $\varphi_h$ is below $10^{-12}$ and that of $u_h$ is below $10^{-10}$.
	The stabilization parameters are $\gamma_u=\gamma_\varphi=1$.
\end{remark}

\subsection{Example 1}
\label{subsec:manufactured-solution}

The computational domain is the unit square $\O=[0,1]^2$, with
non-homogeneous Neumann conditions on $\partial\O$. The exact solution is
\[
u = 1+\tfrac12\sin(\pi x)\cos(\pi y)\cos t,\quad
\p = \cos(\pi x)\sin(\pi y)\sin t,\quad
\o = \bigl(1+\cos(\pi x)\cos(\pi y)\bigr)(t+1).
\]
The source terms in the three equations
of \eqref{MS1}--\eqref{MS3} and the Neumann data $-\fu\cdot\n$, $-\fp\cdot\n$
on $\partial\O$ are computed accordingly. The physical parameters are
$\chi=\mu=\beta=\delta=1$ and $a=2$, and the final time is $T=1$.

The spatial convergence is studied on a sequence of unstructured triangular
meshes generated with \texttt{initmesh} (PDE Toolbox), with maximum edge
length $h\in\{2^{-1},2^{-2},2^{-3},2^{-4}\}$. We fix the  time step 
$dt=2\times10^{-4}$,  so the temporal error is negligible
with respect to the spatial one (%
halving it leaves the errors on
the finest mesh essentially unchanged).

\begin{table}[h!]
	\centering
	\footnotesize
	\setlength{\tabcolsep}{4pt}
	\begin{tabular}{lccccc}
		\toprule
		$h$ & $\|u-u_h\|_\O$ & $\|\p-\p_h\|_\O$ & $\|\o-\o_h\|_\O$ & $\|\fu-\fu_h\|_\O$ & $\|\fp-\fp_h\|_\O$ \\
		\midrule
		\multicolumn{6}{l}{\HDGku} \\
		$2^{-1}$ & 4.748e-2 \qquad - & 8.077e-2 \qquad -& 1.038e-1\qquad - & 8.682e-2\qquad - & 1.505e-1\qquad  -\\
		$2^{-2}$ & 1.532e-2 (1.63) & 2.638e-2 (1.61) & 2.993e-2 (1.79) & 2.847e-2 (1.61) & 4.837e-2 (1.64) \\
		$2^{-3}$ & 3.918e-3 (1.97) & 6.741e-3 (1.97) & 7.592e-3 (1.98) & 7.321e-3 (1.96) & 1.220e-2 (1.99) \\
		$2^{-4}$ & 9.755e-4 (2.01) & 1.685e-3 (2.00) & 1.930e-3 (1.98) & 1.781e-3 (2.04) & 2.990e-3 (2.03) \\
		\midrule
		\multicolumn{6}{l}{\HDGkd} \\
		$2^{-1}$ & 8.154e-3  \qquad -& 1.322e-2  \qquad -& 1.513e-2  \qquad -& 1.711e-2 \qquad - & 2.855e-2 \qquad - \\
		$2^{-2}$ & 1.352e-3 (2.59) & 2.256e-3 (2.55) & 2.246e-3 (2.75) & 2.507e-3 (2.77) & 4.125e-3 (2.79) \\
		$2^{-3}$ & 1.536e-4 (3.14) & 2.640e-4 (3.10) & 2.611e-4 (3.10) & 2.864e-4 (3.13) & 4.805e-4 (3.10) \\
		$2^{-4}$ & 1.825e-5 (3.07) & 3.176e-5 (3.05) & 3.078e-5 (3.08) & 3.384e-5 (3.08) & 5.749e-5 (3.06) \\
		\bottomrule
	\end{tabular}
	\caption{Example 1: $L^2(\O)$-errors at the final time $T=1$ for \HDGku\ and
		\HDGkd, $dt=2\times10^{-4}$. In parentheses, the convergence rate
		between two consecutive meshes.}
	\label{tab:ex1-errors}
\end{table}

Table~\ref{tab:ex1-errors} reports the $L^2(\O)$ errors at $t=T$. All five
quantities converge with the optimal order $\mathcal O(h^{k+1})$ 
in agreement with Theorem~\ref{thm:error}. 

\subsection{Example 2: pattern formation on a square and on a curved domain}
\label{sec:turing-patterns}

\begin{figure}[htbp]
	\centering
	\includegraphics[width=0.8\textwidth]{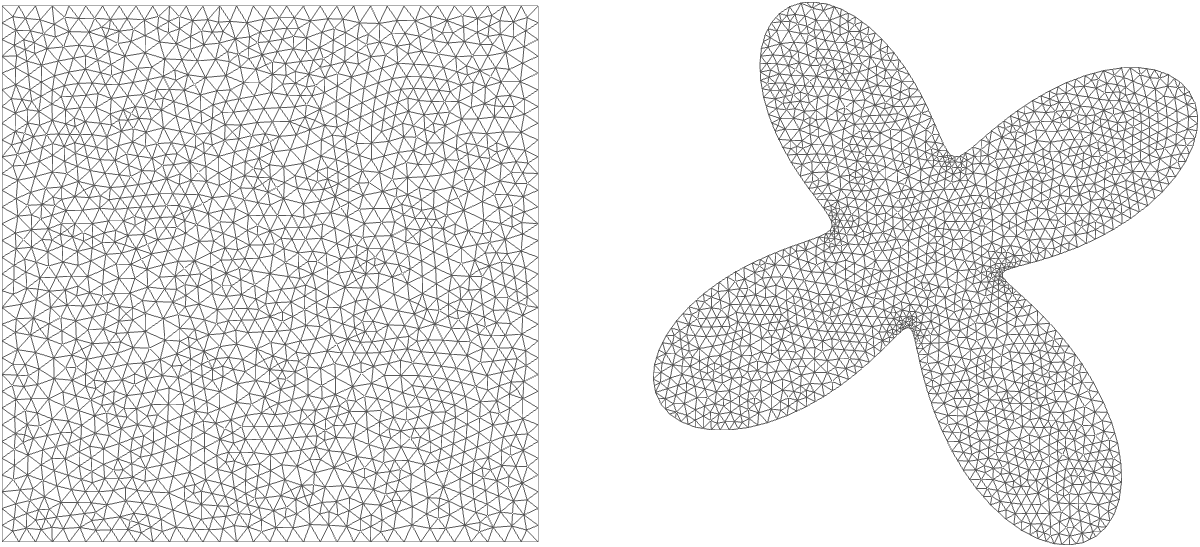}
	\caption{The two domains considered in Example 2 and their meshes with $h\approx0.59$: the square $[0,6\pi]^2$ (left) and the four-leaf domain (right).
	}
	\label{fig:mesh4leaf}
\end{figure}

In this example we show that the method reproduces the Turing patterns of the MS
model computed in \cite{Bisi_pattern} on a square, and that the same computation
carries over, without any change, to a curved, non-convex domain. The two domains
are shown, with the corresponding meshes, in Figure~\ref{fig:mesh4leaf}. The first
one is the square $[0,6\pi]^2$ considered in \cite{Bisi_pattern}. 
The second one, which we refer to as the four-leaf domain, is the region bounded
by the curve $\rho(\omega)=3\pi\bigl(1+0.55\cos(4\omega)+0.1\sin(3\omega)\bigr)$,
$\omega\in[0,2\pi)$, in polar coordinates, rotated by an angle $\pi/6$.
It is a non-convex region with four lobes of different size, whose area, $322$,
is comparable with that of the square ($355$). Its boundary is approximated by
the polygon whose $149$ vertices are obtained by sampling the curve at equispaced
values of $\omega$. Both domains are meshed by the same mesh generator with
$h=6\pi/32\approx0.59$, which gives $3402$ elements for the square and $3556$ for
the four-leaf domain.

With $f$ and $g$ as above and $\beta=\delta=1$, system \eqref{MS1}--\eqref{MS3}
coincides with the dimensionless model of \cite{Bisi_pattern}. The constant
equilibrium $(\bar u,\bar\varphi,\bar{\mathfrak d})=(1,2,1)$ is Turing unstable
for \cite{desvillettes2022global,Bisi_pattern}
\[
\chi>\chi_c=\frac{\bigl(1+\sqrt{\mu(a-1)}\bigr)^2}{\beta f(1)}.
\]
Note that $\chi_c$ is larger than the threshold $\chi_{subcrit}$ of
Theorem~\ref{thm:stability}, which is only a sufficient condition for stability;
the two coincide for $\mu(a-1)=1$ \cite[Remark~1.4]{desvillettes2022global}.

We consider the three parameter sets of \cite[Fig.~3]{Bisi_pattern}, which lead to
stripes, hexagons and spots, respectively:
\begin{center}
	\begin{tabular}{lccc}
		\hline
		pattern & $\mu$ & $\chi$ & $\chi_c$\\
		\hline
		stripes  & 0.8  & 8.33 & 7.18\\
		hexagons & 0.8  & 7.32 & 7.18\\
		spots    & 0.02 & 2.62 & 2.61\\
		\hline
	\end{tabular}
\end{center}

All runs use HDG$(2,2)$ on the meshes of Figure~\ref{fig:mesh4leaf}, homogeneous Neumann conditions, 
BDF2 with $\Delta t=0.1$ and the Newton linearization
(Remark~\ref{rem:implementation}). The initial datum is the equilibrium plus, in
each component, an independent realization of a smooth random perturbation 
$\varepsilon\,r(x,y)$,
	with
\[
r(x,y)=\sum_{j=1}^{50}a_j\cos(\mathbf k_j\cdot\mathbf x+\theta_j),
\]
where the wave vectors $\mathbf k_j\in[-k_{\max},k_{\max}]^2$, the phases
$\theta_j$ and the amplitudes $a_j\ge0$ are drawn
with a fixed seed, and
$\sum_j a_j=1/2$, so that $|r|\le 1/2$; we take $\varepsilon=10^{-3}$ and
$k_{\max}=2$ for stripes and hexagons, and $\varepsilon=10^{-2}$ and $k_{\max}=4$
for spots, whose critical wavenumber is larger and linear growth rate smaller.
Each run is carried out in blocks of 500 time units and is stopped at
the end of the first block at which the relative change
$\|u_h(t)-u_h(t-10)\|_\Omega/\|u_h(t)\|_\Omega$ is below $10^{-4}$ (some runs were
continued further, to check stationarity), and in any case at $t=3000$.

\

\begin{figure}[htbp]
	\centering
	\includegraphics[width=0.8\textwidth]{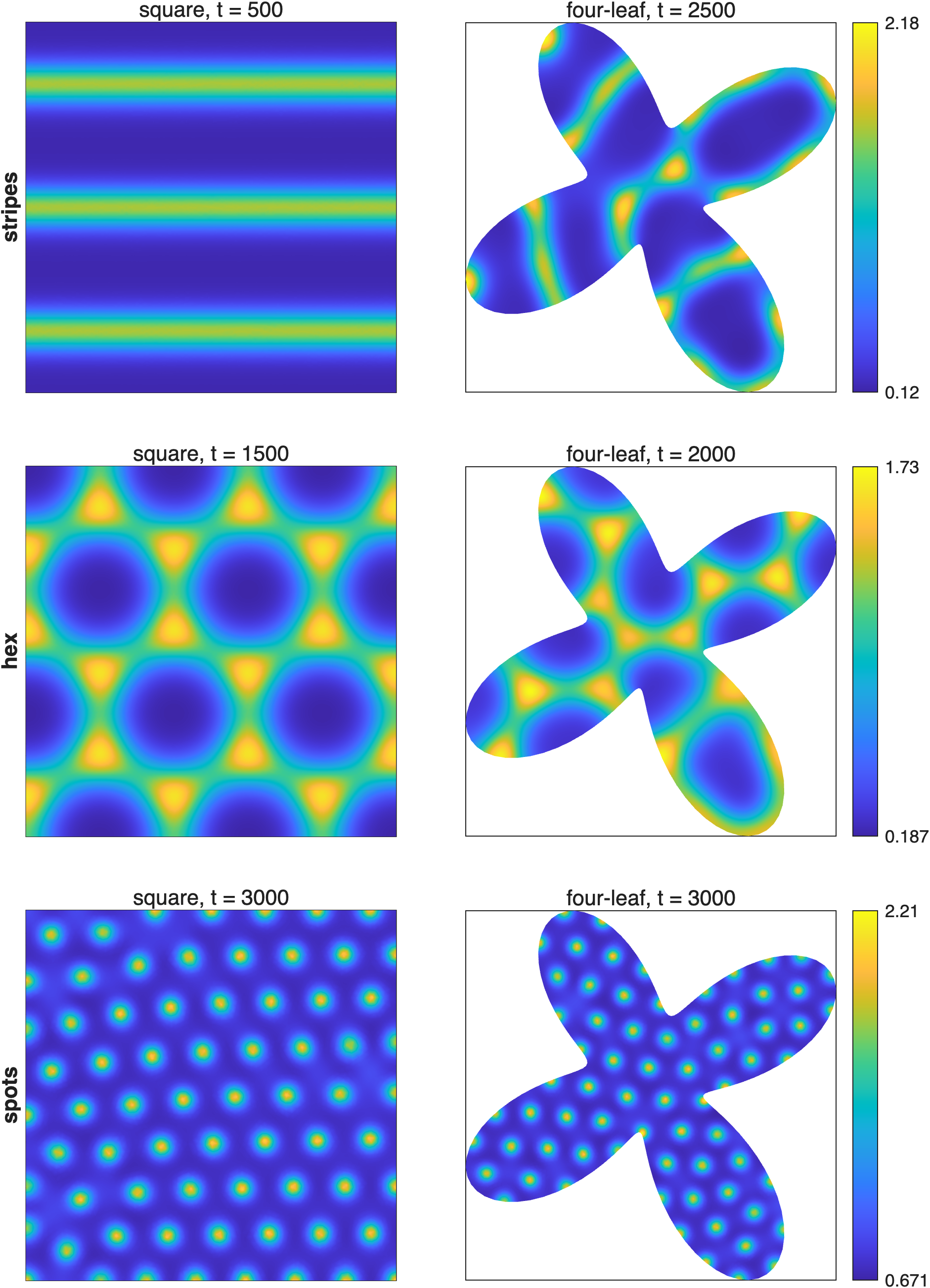}
	\caption{
		Example 2:  $u_h$ at the final time of Table~\ref{tab:ex2-final} for the
		stripes (top), hexagon (middle) and spot (bottom) parameters of
		\cite[Fig.~3]{Bisi_pattern}, on the square $[0,6\pi]^2$ (left) and on the
		four-leaf domain (right; not drawn to the same scale as the square,
		its bounding box being about 40\% larger).
		 The color scale is common to the two panels of a
		row.
		}
	\label{fig:ex2-final}
\end{figure}

\paragraph{\bf Square domain.}
On the square the three patterns of \cite[Fig.~3]{Bisi_pattern} are recovered
(Figure~\ref{fig:ex2-final}, left column), with density ranges comparable to the
ones reported there.
The parameters select the type
and the wavelength of the pattern, but not its orientation or the position of its
defects, which depend on the initial perturbation. In particular, for the hexagon
parameters the random perturbation alone leads to a lattice with one defect; 
the regular lattice of \cite[Fig.~3(b)]{Bisi_pattern} shown here is obtained by adding
to $u_0$ a perturbation of amplitude $10^{-2}$ along the three wave vectors
$(0,6)\pi/L$, $(5,\pm3)\pi/L$, $L=6\pi$, compatible with the no-flux conditions, the first mode having sign opposite to
the other two.

\

\paragraph{\bf A curved domain.}
The right column of Figure~\ref{fig:ex2-final} shows the same three runs on the four-leaf domain. The initial perturbation is the random field alone.
Nothing changes in the method: the no-flux conditions are
imposed on the boundary edges through the numerical fluxes, exactly as the
transmission conditions on the interior edges, and the geometry enters only
through the mesh. The same three patterns form, with comparable wavelengths, and
they adapt to the geometry: the stripes meet the boundary at a right angle or run
along it, the hexagonal network follows the lobes, and the spots form a nearly
hexagonal lattice with defects where the lobes meet.

Table~\ref{tab:ex2-final} reports, for the six runs, the final time, the relative
change of $u_h$ over the last 10 time units, the range of $u_h$ and the average
number of Newton iterations per step. Stripes and hexagons reach a stationary
state on both domains. For the spots the relative change levels off at
$2$--$3\cdot10^{-4}$, due to the slow motion of the lattice defects, and the runs
were stopped at $t=3000$.

\begin{table}[htbp]
	\centering
	\footnotesize
	\setlength{\tabcolsep}{4pt}
	\begin{tabular}{llrccc}
			\toprule
			pattern & domain & $T$ & rel.\ change at $T$ & $u_h(T)$ & Newton it./step \\
			\midrule
			stripes  & square    &  500 & $3.4\cdot10^{-5}$  & $[0.140, 1.610]$ & 3.0      \\
			stripes  & four-leaf & 2500 & $2.1\cdot10^{-5}$  & $[0.120, 2.180]$ & 2.0--3.2  \\
			hexagons & square    & 1500 & $9.5\cdot10^{-15}$ & $[0.187, 1.636]$ & 1.7     \\
			hexagons & four-leaf & 2000 & $3.5\cdot10^{-12}$ & $[0.232, 1.731]$ & 2.4   \\
			spots    & square    & 3000 & $2.2\cdot10^{-4}$  & $[0.696, 2.089]$ & 2.0--3.0  \\
			spots    & four-leaf & 3000 & $2.9\cdot10^{-4}$  & $[0.671, 2.210]$ & 2.8--3.0  \\
			\bottomrule
		\end{tabular}
	\caption{Example 2: final time, stationarity, range of $u_h$ and
			average
			number of Newton iterations per step  
			 for the six runs of Figure~\ref{fig:ex2-final} (\HDGkd, $h=6\pi/32$,
			$dt=0.1$). The relative change is
			$\|u_h(T)-u_h(T-10)\|_\O/\|u_h(T)\|_\O$. }
	\label{tab:ex2-final}
\end{table}

\section{Conclusions}\label{sec:conclusions}

We have presented an equal-order hybridizable discontinuous Galerkin method for a
chemotaxis model of multiple sclerosis, in which a Keller--Segel type system is
coupled with an ordinary differential equation for the tissue damage. The damage
variable is treated element by element and adds no globally coupled unknown, so
that the method has the same condensed structure as the HDG method for the
Keller--Segel system of \cite{BBPP}. For the semi-discrete problem we proved
existence and uniqueness of the discrete solution without assuming its positivity,
a stability bound uniform in the mesh size for bounded chemotactic sensitivity,
and an a priori error estimate which is of optimal order $k+1$ for the densities
and for the fluxes when the solution is sufficiently smooth. The numerical tests
confirm these rates for $k=1,2$.

The pattern formation experiments show that the method reproduces the 
patterns known for the model on a square and that, since the geometry and the
boundary conditions enter only through the mesh and the numerical fluxes, the same
computation applies unchanged to a curved, non-convex domain. This makes the
method a suitable tool for the simulation of this class of models on domains of
general shape, such as those arising from the anatomy of the brain.

Several issues are left to future work. The analysis presented here concerns the
semi-discrete problem, while the computations use
implicit time discretizations
(backward Euler and BDF2) with Picard or Newton linearizations of the nonlinear
terms, whose analysis is still to be
carried out. The discrete solution is not guaranteed to be non-negative;
in the
pattern formation experiments $u_h$ was positive at the final time of all runs
(Table~\ref{tab:ex2-final}),
but positivity preserving variants deserve investigation. Curved boundaries are
approximated by polygons, which is adequate for the purposes of this paper but
limits the accuracy for high polynomial degrees. Finally, the extension of the
method to polytopal meshes, and its application to realistic three dimensional
geometries, will be the subject of a forthcoming work.

\bibliographystyle{plain}
\bibliography{paper_arxiv}

\end{document}